\documentclass[11pt]{article}
\usepackage{amsmath}

\usepackage[margin=1in, centering]{geometry}

\usepackage{float}

\usepackage{placeins}
\usepackage{verbatim}

\usepackage{authblk}
\usepackage{url}
\usepackage{amsthm}

\usepackage[dvips]{graphicx}
\usepackage{color}
\usepackage{amsfonts,amssymb,amsmath,}
\usepackage{mathrsfs}
\usepackage[T1]{fontenc}
\usepackage{ae}

\newtheorem{theorem}{Theorem}

\newtheorem{coro}{Corollary}[theorem]

\newtheorem{lem}{Lemma}

\numberwithin{equation}{section}
\numberwithin{table}{section}
\numberwithin{figure}{section}

\renewcommand{\(}{\left(}
\renewcommand{\)}{\right)}

\allowdisplaybreaks[4]

\begin{document}
\title{Chebyshev's bias for products of $k$ primes}
\author{Xianchang Meng}
\date{}

\maketitle
\begin{abstract}
For any $k\geq 1$, we study the distribution of the difference between the number of integers $n\leq x$ with $\omega(n)=k$  or $\Omega(n)=k$ in two different arithmetic progressions, where $\omega(n)$ is the number of distinct prime factors of $n$ and $\Omega(n)$ is the number of prime factors of $n$ counted with multiplicity . Under some reasonable assumptions, we show that, if $k$ is odd, the integers with $\Omega(n)=k$ have preference for quadratic non-residue classes; and if $k$ is even, such integers have preference for quadratic residue classes. This result confirms a conjecture of Richard Hudson. However, the integers with $\omega(n)=k$ always have preference for quadratic residue classes. Moreover, as $k$ increases, the biases become smaller and smaller for both of the two cases.


\end{abstract}

\let\thefootnote\relax

\let\thefootnote\relax

\section{Introduction and statement of results}
First, we consider products of $k$ primes in arithmetic progressions. Let
\begin{equation*}
\pi_k(x; q, a)=|\{n\leq x: \omega(n)=k, ~n \equiv a \bmod q\}|,
\end{equation*}
and
\begin{equation*}
N_k(x; q, a)=|\{n\leq x: \Omega(n)=k, ~n \equiv a \bmod q\}|,
\end{equation*}
 where $\omega(n)$ is the number of distinct prime divisors of $n$, and $\Omega(n)$ is the number of prime divisors of $n$ counted with multiplicity.
For example, when $k=1$, $N_1(x; q, a)$ is the number of primes $\pi(x; q, a)$ in the arithmetic progression $a \bmod q$; and $\pi_1(x; q, a)$ counts the number of prime powers $p^l\leq x$ for all $l\geq 1$ in the arithmetic progression $a \bmod q$.

Dirichlet (1837) \cite{Diri}  showed that, for any $a$ and $q$ with $(a, q)=1$, there are infinitely many primes in the arithmetic progression $a \bmod q$. Moreover, for any $(a, q)=1$,
\begin{equation*}
\pi(x; q, a)\sim \frac{x}{\phi(q) \log x},
\end{equation*}
where $\phi$ is Euler's totient function \cite{Dave}.
Analogous asymptotic formulas are available for products of $k$ primes. Landau (1909) \cite{Lan} showed that, for each fixed integer $k\geq 1$,
\begin{equation*}
	N_k(x):=\left| \{n\leq x: \Omega(n)=k \}  \right|\sim \frac{x}{\log x}\frac{(\log\log x)^{k-1}}{(k-1)!}.
\end{equation*}
The same asymptotic is also true for the function $\pi_k(x):=| \{n\leq x: \omega(n)=k \}  |$.  For more precise formulas, see \cite{Tene} (II. 6, Theorems 4 and 5). Using similar methods as in \cite{Dave} and \cite{Tene}, one can show that, for any fixed residue class $a \bmod q$ with $(a, q)=1$,
\begin{equation*}
	N_k(x; q, a)\sim \pi_k(x; q, a)\sim \frac{1}{\phi(q)} \frac{x}{\log x} \frac{(\log\log x)^{k-1}}{(k-1)!}.
\end{equation*}

 For the case of counting primes ($\Omega(n)=1$), Chebyshev (1853) \cite{Che} observed that there seem to be more primes in the progression $3\bmod 4$ than in the progression $1\bmod 4$. That is, it appears that $\pi(x; 4,3)\geq \pi(x; 4, 1)$.
 In general, for any $a\not\equiv b \bmod q$ and $(a, q)=(b, q)=1$, one can study the behavior of the functions
 \begin{align*}
 \Delta_{\omega_k}(x; q, a, b)&:=\pi_k(x; q, a)-\pi_k(x; q, b),\\
 \Delta_{\Omega_k}(x; q, a, b)&:=N_k(x; q, a)-N_k(x; q, b).
 \end{align*}
Denote $\Delta(x; q, a, b):=\Delta_{\Omega_1}(x; q, a, b)$.
  Littlewood \cite{Lit} proved that $\Delta(x; 4, 3, 1)$ changes sign infinitely often. Actually, $\Delta(x; 4, 3, 1)$ is negative for the first time at $x=26, 861$ \cite{Lee}. Knapowski and Tur{\'a}n published a series of papers starting with \cite{KT} about the sign changes and extreme values of the functions $\Delta(x; q, a, b)$. And such problems are colloquially known today as "prime race problems". Irregularities in the distribution, that is, a tendency for $\Delta(x; q, a, b)$ to be of one sign is known as "Chebyshev's bias". For a nice survey of such works, see \cite{F-K} and \cite{Gran-Mar}.

Chebyshev's bias can be well understood in the sense of \textit{logarithmic density}. We say a set $S$ of positive integers has \textit{logarithmic density}, if the following limit exists:
\begin{equation*}
\delta(S)=\lim_{x\rightarrow\infty} \frac{1}{\log x} \sum_{\substack{n\leq x\\ n\in S}} \frac{1}{n}.
\end{equation*}
Let $\delta_{f_k}(q; a, b)=\delta(P_{f_k}(q; a, b))$, where $P_{f_k}(q; a, b)$ is the set of integers with $\Delta_{f_k}(n; q, a, b)>0$, and $f=\Omega {~\rm or~} \omega$.
In order to study the Chebyshev's bias and the existence of the logarithmic density, we need the following assumptions:

	1) the Extended Riemann Hypothesis (${\rm ERH_q}$) for Dirichlet L-functions modulo $q$;
	
	2) the Linear Independence conjecture (${\rm LI_q}$), the imaginary parts of the zeros of all Dirichlet L-functions modulo $q$ are linearly independent over $\mathbb{Q}$.

Under these two assumptions, Rubinstein and Sarnak \cite{RS} showed that, for Chebyshev's bias for primes ($\Omega(n)=1$), the logarithmic density $\delta_{\Omega_1}(q;a, b)$ exists, and in particular, $\delta_{\Omega_1}(4; 3, 1)\approx 0.996$ which indicates a strong bias for primes in the arithmetic progression $3 \bmod 4$.
Recently, using the same assumptions, Ford and Sneed \cite{Ford-S} studied the Chebyshev's bias for products of two primes with $\Omega(n)=2$ by transforming this problem into manipulations of some double integrals.  They connected $\Delta_{\Omega_2}(x; q, a, b)$ with $\Delta(x; q, a, b)$, and showed that $\delta_{\Omega_2}(q; a, b)$ exists and the bias is in the opposite direction to the case of primes, in particular, $\delta_{\Omega_2}(4; 3, 1)\approx 0.10572$ which indicates a strong bias for the arithmetic progression $1 \bmod 4$. 

By orthogonality of Dirichlet characters, we have
\begin{equation}\label{Delta-Omega}
\Delta_{\Omega_k}(x; q, a, b)=\frac{1}{\phi(q)} \sum_{\chi\neq \chi_0 \bmod q} (\overline{\chi}(a)-\overline{\chi}(b))\sum_{\substack{n\leq x\\ \Omega(n)=k}}\chi(n),
\end{equation}
and
\begin{equation}\label{Delta-omega}
\Delta_{\omega_k}(x; q, a, b)=\frac{1}{\phi(q)} \sum_{\chi\neq \chi_0 \bmod q} (\overline{\chi}(a)-\overline{\chi}(b))\sum_{\substack{n\leq x\\ \omega(n)=k}}\chi(n).
\end{equation}
The inner sums over $n$ are usually analyzed using analytic methods. Neither the method of Rubinstein and Sarnak \cite{RS} nor the method of Ford and 
Sneed \cite{Ford-S} readily generalizes to handle the cases of more prime factors 
($k\geq 3$).  From the point of view of $L$-functions, the  most 
natural sum to consider is 
\begin{equation}\label{Lambda-analog}
\sum_{\substack{n_1\cdots n_k\leq x\\n_1\cdots n_k \equiv a \bmod q}} \Lambda(n_1)\cdots\Lambda(n_k).
\end{equation}
However, estimates for $\Delta_{\Omega_k}(x; q, a, b)$ or $\Delta_{\omega_k}(x; q, a, b)$ cannot be readily recovered from such an analogue by partial summation. 
Ford and Sneed \cite{Ford-S} overcome this obstacle in the case $k=2$ by means of the 2-dimensional integral 
$$ 
\int_0^\infty \int_0^\infty \sum_{p_1 p_2\le x} \frac{\chi(p_1 p_2)\log p_1 \log p_2}{p_1^{u_1} p_2^{u_2}} du_1 du_2. 
$$ 
Analysis of an analogous $k$-dimensional integral leads to an explosion of 
cases, depending on the relative sizes of the variables $u_j$, and becomes 
increasingly messy as $k$ increases.


We take an entirely different approach, working directly with the unweighted 
sums. 
We express the associated Dirichlet series in terms of products of the logarithms of Dirichlet $L$-functions, then apply Perron's formula, and use Hankel contours to avoid the zeros of $L(s, \chi)$ and the point $s=\frac{1}{2}$.
Using the same assumptions 1) and 2), we show that, for any $k\geq 1$, both $\delta_{\Omega_k}(q; a, b)$ and $\delta_{\omega_k}(q; a, b)$ exist.  Moreover, we show that, as $k$ increases, if $a$ is a quadratic non-residue and $b$ is a quadratic residue, the bias oscillates with respect to the parity of $k$ for the case $\Omega(n)=k$, but $\delta_{\omega_k}(q; a, b)$ increases from below $\frac{1}{2}$ monotonically.

For some of our results, we need only a much weaker substitute for condition ${\rm LI_q}$, which we call the Simplicity Hypothesis (${\rm SH_q}$):
$\forall \chi\neq \chi_0 \bmod q$, $L(\frac{1}{2}, \chi)\neq 0$ and the zeros of $L(s, \chi)$ are simple.
Let 
	$$N(q, a):=\#\{ u \bmod q: u^2\equiv a \bmod q\}.$$
Then, using the weaker assumptions ${\rm SH_q}$ and ${\rm ERH_q}$, we prove the following theorems.
\begin{theorem}\label{thm-Delta-Omega}
	Assume $ERH_q$ and ${SH_q}$. Then, for any fixed $k\geq 1$, and fixed large $T_0$,
	\begin{align*}
	\Delta_{\Omega_k}(x; q, a, b)&=\frac{1}{(k-1)!} \frac{\sqrt{x} (\log\log x)^{k-1}}{\log x}\Bigg\{ \frac{(-1)^k}{\phi(q)}\sum_{\chi\neq \chi_0}\(\overline{\chi}(a)-\overline{\chi}(b)\)\sum_{\substack{|\gamma_{\chi}|\leq T_0\\L(\frac{1}{2}+i\gamma_{\chi}, \chi)=0}} \frac{x^{i\gamma_{\chi}}}{\frac{1}{2}+i\gamma_{\chi}}\nonumber\\
	&\quad +\frac{(-1)^k}{2^{k-1}}\frac{N(q,a)-N(q,b)}{\phi(q)} +\Sigma_k(x; q, a, b, T_0)\Bigg\},
	\end{align*}
	where 
	\begin{equation}\nonumber
		\limsup_{Y\rightarrow\infty}\frac{1}{Y} \int_1^Y \left|\Sigma_k(e^y; q, a, b, T_0) \right|^2 dy\ll \frac{\log^2 T_0}{T_0}.
	\end{equation}
	\end{theorem}
Since $\Delta_{\Omega_1}(x; q, a, b)=\Delta(x; q, a, b)$, we get the following corollary. 	
	\begin{coro}\label{Coro-Omega}
		Assume $ERH_q$ and ${SH_q}$. Then, for any fixed $k\geq 2$,
		\begin{align*}
		\frac{\Delta_{\Omega_k}(x; q, a, b)\log x}{\sqrt{x} (\log\log x)^{k-1}}&=\frac{(-1)^{k+1}}{(k-1)!} \(1-\frac{1}{2^{k-1}} \) \frac{N(q,a)-N(q, b)}{ \phi(q)}\nonumber\\
		&\quad +\frac{(-1)^{k+1}}{(k-1)!}\frac{\Delta(x; q, a, b) \log x}{ \sqrt{x}} +\Sigma'_k(x; q, a, b),	
		\end{align*}
		where, as $Y\rightarrow \infty$,
		\begin{equation}\nonumber
		\frac{1}{Y} \int_1^Y |\Sigma'_k(e^y; q, a, b) |^2 dy=o(1).
		\end{equation}
	\end{coro}

\begin{figure}[htp]
	\centering
	\includegraphics[width=15cm, height=10cm]{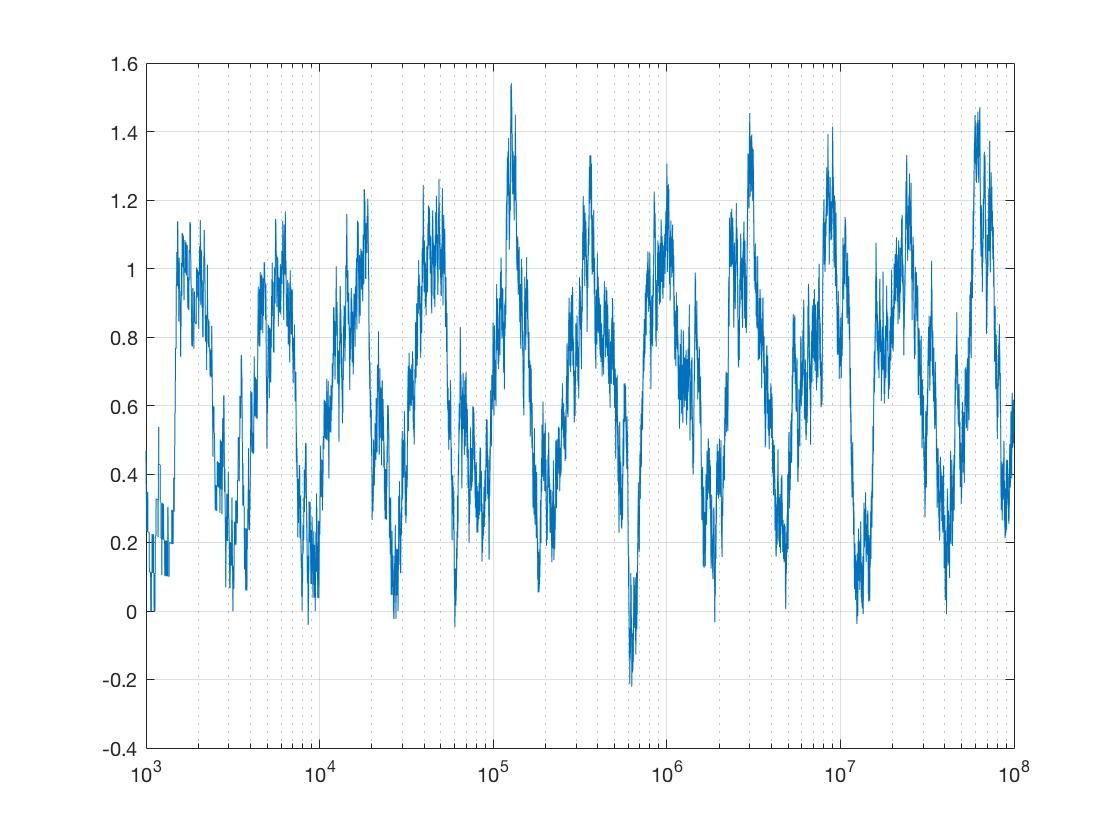}\\
	\caption{$\frac{2\log x}{\sqrt{x}(\log\log x)^2}\Delta_{\Omega_3}(x; 4, 3, 1)$}\label{fig-Big-Omega-3}
\end{figure}
\begin{figure}[htp]
	\centering
	\includegraphics[width=15cm, height=10cm]{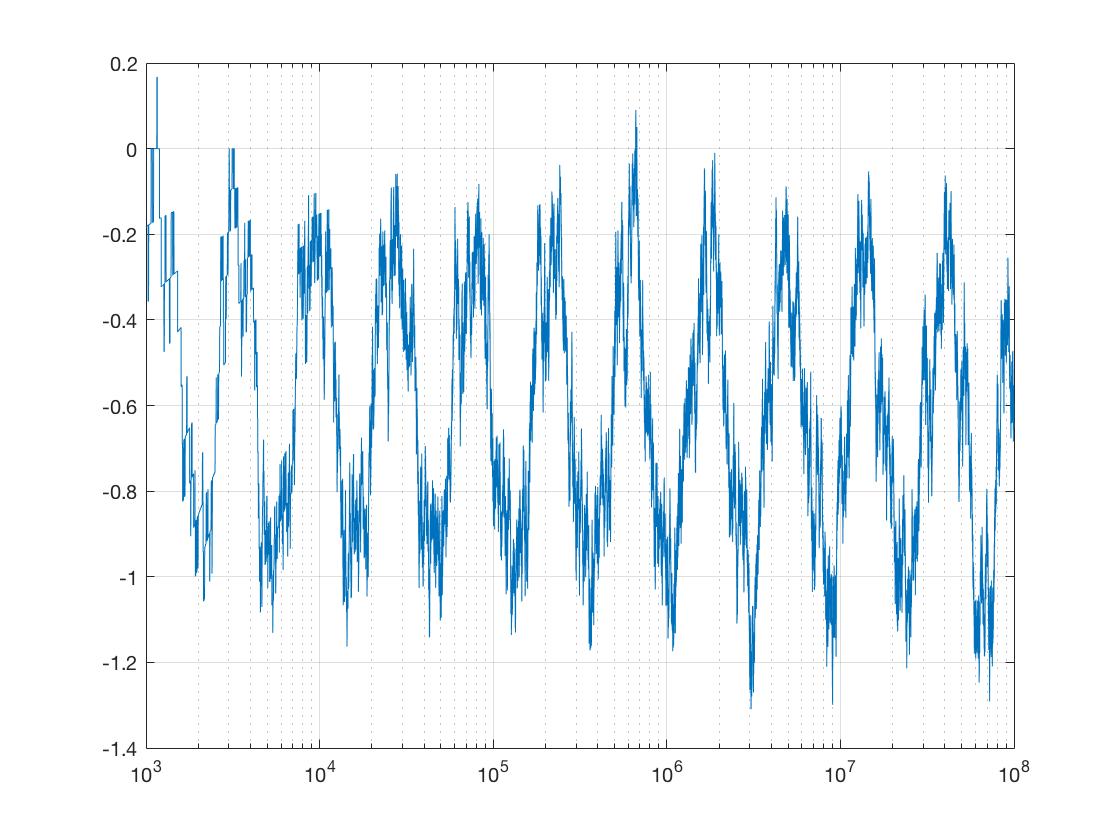}\\
	\caption{$\frac{6\log x}{\sqrt{x}(\log\log x)^3}\Delta_{\Omega_4}(x; 4, 3, 1)$}\label{fig-Big-Omega-4}
\end{figure}

In the above theorem, the constant $\frac{(-1)^k}{2^{k-1}}\frac{N(q,a)-N(q,b)}{\phi(q)}$ represents the bias in the distribution of products of $k$ primes counted with multiplicity. Richard Hudson conjectured that, as $k$ increases, the bias would change directions according to the parity of $k$. Our result above confirms his conjecture (under $\rm ERH_q$ and $\rm SH_q$).  Figures \ref{fig-Big-Omega-3} and \ref{fig-Big-Omega-4} show the graphs corresponding to $(q, a, b)=(4, 3, 1)$ for  $\frac{2\log x}{\sqrt{x}(\log\log x)^2}\Delta_{\Omega_3}(x; 4, 3, 1)$ and $\frac{6\log x}{\sqrt{x} (\log\log x)^3}\Delta_{\Omega_4}(x; 4, 3, 1)$, plotted on a logarithmic scale from $x=10^3$ to $x=10^8$.  In these graphs, the functions do not appear to be oscillating around $\frac{1}{4}$ and $-\frac{1}{8}$ respectively as predicted in our theorem. This is caused by some terms of order $\frac{1}{\log\log x}$ and even lower order terms, and $\log\log 10^8\approx 2.91347$ and $\frac{1}{\log\log 10^8}\approx 0.343233$. However, we can still observe the expected direction of the bias through these graphs.

\bigskip
For the distribution of products of $k$ primes counted without multiplicity, we have the following theorem. In this case, the bias will be determined by the constant $\frac{N(q,a)-N(q,b)}{2^{k-1}\phi(q)}$ in the theorem below. 
\begin{theorem}\label{thm-Delta-omega}
	Assume $ERH_q$ and ${SH_q}$. Then, for any fixed $k\geq 1$, and fixed large $T_0$,
	\begin{align*}
	\Delta_{\omega_k}(x; q, a, b)&=\frac{1}{(k-1)!} \frac{\sqrt{x} (\log\log x)^{k-1}}{\log x}\Bigg\{ \frac{(-1)^k}{\phi(q)}\sum_{\chi\neq \chi_0}\(\overline{\chi}(a)-\overline{\chi}(b)\)\sum_{\substack{|\gamma_{\chi}|\leq T_0\\L(\frac{1}{2}+i\gamma_{\chi}, \chi)=0}} \frac{x^{i\gamma_{\chi}}}{\frac{1}{2}+i\gamma_{\chi}}\nonumber\\
	&\quad +\frac{N(q,a)-N(q,b)}{2^{k-1}\phi(q)} +\widetilde{\Sigma}_k(x; q, a, b, T_0)\Bigg\},
	\end{align*}
	where 
	\begin{equation}\nonumber
		 \limsup_{Y\rightarrow\infty} \frac{1}{Y} \int_1^Y \left|\widetilde{\Sigma}_k(e^y; q, a, b, T_0) \right|^2 dy\ll \frac{\log^2 T_0}{T_0}.
	\end{equation}
\end{theorem}

	\begin{coro}\label{Coro-omega}
		Assume $ERH_q$ and ${SH_q}$. Then, for any fixed $k\geq 1$,
		\begin{align*}
		\frac{\Delta_{\omega_k}(x; q, a, b)\log x}{\sqrt{x} (\log\log x)^{k-1}}&= \(\frac{1}{2^{k-1}}+(-1)^{k+1} \) \frac{N(q,a)-N(q, b)}{(k-1)! \phi(q)}\nonumber\\
		&\quad +\frac{(-1)^{k+1}}{(k-1)!}\frac{\Delta(x; q, a, b) \log x}{ \sqrt{x}} +\widetilde{\Sigma}'_k(x; q, a, b),	
		\end{align*}
		where, as $Y\rightarrow \infty$,
		\begin{equation}\nonumber
		\frac{1}{Y} \int_1^Y |\widetilde{\Sigma}'_k(e^y; q, a, b) |^2 dy=o(1).
		\end{equation}
	\end{coro}


For the distribution of $\Delta(x; q, a, b)$, Rubinstein and Sarnak \cite{RS} showed the following theorem. This is the version from \cite{Ford-S}.
\newtheorem*{theoremRS}{Theorem RS}
\begin{theoremRS}\label{thm-RS}
	Assume $ERH_q$ and $LI_q$. For any $a\not\equiv b \bmod q$ and $(a, q)=(b, q)=1$, the function
	\begin{equation}\nonumber
	\frac{u\Delta(e^u; q, a, b)}{e^{u/2}}
	\end{equation}
	has a probabilistic distribution. This distribution i) has mean $\frac{N(q, b)-N(q, a)}{\phi(q)}$, ii) is symmetric with respect to its mean, and iii) has a continuous density function.
\end{theoremRS}

Corollaries \ref{Coro-Omega}, \ref{Coro-omega}, and Theorem RS imply the following result.

\begin{theorem}\label{thm-density-exist}
	Let $a\not\equiv b \bmod q$ and $(a, q)=(b, q)=1$. Assuming $ERH_q$ and $LI_q$, for any $k\geq 1$, $\delta_{\Omega_k}(q; a, b)$ and $\delta_{\omega_k}(q; a, b)$ exist. More precisely, if $a$ and $b$ are both quadratic residues or both quadratic non-residues, then $\delta_{\Omega_k}(q; a, b)=\delta_{\omega_k}(q; a, b)=\frac{1}{2}$. Moreover,
if $a$ is a quadratic non-residue and $b$ is a quadratic residue, then, for any $k\geq 1$,
	$$1-\delta_{\Omega_{2k-1}}(q; a, b)<\delta_{\Omega_{2k}}(q; a, b)< \frac{1}{2}<\delta_{\Omega_{2k+1}}(q; a, b)<1-\delta_{\Omega_{2k}}(q; a, b),$$
    $$\delta_{\omega_{k}}(q; a, b)<\delta_{\omega_{k+1}}(q; a, b)<\frac{1}{2},$$
	$$\delta_{\Omega_{2k}}(q; a, b)=\delta_{\omega_{2k}}(q; a, b), \qquad \delta_{\Omega_{2k-1}}(q; a, b)+\delta_{\omega_{2k-1}}(q; a, b)=1.$$	
\end{theorem}
\noindent \textbf{Remark 1.} The above results confirm a conjecture of Richard Hudson proposed years ago in his communications with Ford. Borrowing the methods from \cite{RS} (Section 4), we are able to calculate $\delta_{\Omega_k}(q; a, b)$ and $\delta_{\omega_k}(q; a, b)$ precisely for special values of $q$, $a$, and $b$. In particular, we record in Tables \ref{LD-L3} and \ref{LD-L4} the logarithmic densities up to products of 10 primes for two cases: $q=3$, $a=2$, $b=1$, and $q=4$, $a=3$, $b=1$.

\begin{table}[h!]
	\parbox[c]{0.49\textwidth}{\centering
		\begin{tabular}{|c| c|c|}
			\hline
			\multicolumn{3}{|c|}{$q=3$, $a=2$, $b=1$} \\ \hline
			$k$ & $\delta_{\Omega_{k}}(3; 2, 1)$ & $\delta_{\omega_{k}}(3; 2, 1)$\\ \hline
			1  & 0.99906, \cite{RS} & 0.00094 \\ \hline
			2& 0.069629 & 0.069629 \\ \hline
			3 & 0.766925 &  0.233075 \\ \hline
			4 & 0.35829   & 0.35829\\ \hline
			5 &  0.571953  &  0.428047 \\ \hline
			6 & 0.463884   & 0.463884\\ \hline
			7 & 0.518075   &  0.481925\\ \hline
			8 & 0.49096   & 0.49096\\ \hline
			9 & 0.50452   &  0.49548\\ \hline
			10 & 0.49774   & 0.49774\\ \hline
			
		\end{tabular}\caption{$\delta_{\Omega_{k}}(3; 2, 1)$ and $\delta_{\omega_{k}}(3; 2, 1)$}\label{LD-L3}
		}
	\hfill
		\parbox[c]{0.49\textwidth}{\centering
			\begin{tabular}{|c| c|c|}
				\hline
				\multicolumn{3}{|c|}{$q=4$, $a=3$, $b=1$}  \\ \hline
				$k$ & $\delta_{\Omega_{k}}(4; 3, 1)$ & $\delta_{\omega_{k}}(4; 3, 1)$\\ \hline
				1	&  0.9959, \cite{RS}&  0.0041  \\ \hline
				2	&  0.10572, \cite{Ford-S} & 0.10572 \\ \hline
				3 & 0.730311  & 0.269689 \\ \hline
				4 & 0.380029  & 0.380029 \\ \hline
				5 &  0.56061   &  0.43939 \\ \hline
				6 & 0.469616   & 0.469616 \\ \hline
				7& 0.515202   &  0.484798 \\ \hline
				8 & 0.492398  &  0.492398 \\ \hline
				9 & 0.503801   &  0.496199 \\ \hline
				10 & 0.498099  & 0.498099 \\ \hline
			\end{tabular}
			\caption{$\delta_{\Omega_{k}}(4; 3, 1)$ and $\delta_{\omega_{k}}(4; 3, 1)$}\label{LD-L4}   }
\end{table}

For fixed $q$ and large $k$, we give asymptotic formulas for $\delta_{\Omega_k}(q; a, b)$ and $\delta_{\omega_k}(q; a, b)$.
\begin{theorem}\label{thm-density-asym}
	Assume $ERH_q$ and $LI_q$. Let $A(q)$ be the number of real characters $\bmod~q$. Let $a$ be a quadratic non-residue and $b$ be a quadratic residue, and $(a, q)=(b, q)=1$. Then, for any nonnegative integer $K$, and any $\epsilon>0$,
	\begin{equation}\label{formu-delta-Omega}
	\delta_{\Omega_k}(q; a, b)=\frac{1}{2}+\frac{(-1)^{k-1}}{2\pi}\sum_{j=0}^K \(\frac{1}{2^{k-1}} \)^{2j+1}\frac{(-1)^j  A(q)^{2j+1}C_j(q; a,b)}{(2j+1)!} +O_{q, K, \epsilon}\(\frac{1}{(2^{k-1})^{2K+3-\epsilon}}\), 
	\end{equation}
	\begin{equation}\label{formu-delta-omega}
	\delta_{\omega_k}(q; a, b)=\frac{1}{2}-\frac{1}{2\pi}\sum_{j=0}^K \(\frac{1}{2^{k-1}} \)^{2j+1}\frac{(-1)^j  A(q)^{2j+1}C_j(q; a,b)}{(2j+1)!} +O_{q, K, \epsilon}\(\frac{1}{(2^{k-1})^{2K+3-\epsilon}}\), 
	\end{equation}
	where $C_j(q; a, b)$ is some constant depending on $j$, $q$, $a$, and $b$.
	In particular, for $K=0$,
	\begin{equation*}
		\delta_{\Omega_k}(q; a, b)=\frac{1}{2}+(-1)^{k-1}\frac{A(q)C_0(q; a, b)}{2^k\pi} +O_{q, \epsilon}\(\frac{1}{(2^{k})^{3-\epsilon}}\),
		\end{equation*}
		\begin{equation*}
		\delta_{\omega_k}(q; a, b)=\frac{1}{2}-\frac{A(q)C_0(q; a, b)}{2^k\pi} +O_{q, \epsilon}\(\frac{1}{(2^{k})^{3-\epsilon}}\).
	\end{equation*}

\end{theorem}

\noindent \textbf{Remark 2.} 
We have a formula for $C_j(q; a, b)$, 
$$C_j(q; a, b)=\int_{-\infty}^{\infty} x^{2j} \Phi_{q; a, b}(x)dx,$$
where $$\Phi_{q; a, b}(z)=\prod_{\chi\neq \chi_0} \prod_{\substack{\gamma_{\chi}>0\\L(\frac{1}{2}+i\gamma_{\chi})=0}}J_0\(\frac{2|\chi(a)-\chi(b)|z}{\sqrt{\frac{1}{4}+\gamma_{\chi}^2}} \),$$
and $J_0(z)$ is the Bessel function, 
$$J_0(z)=\sum_{m=0}^{\infty} \frac{(-1)^m (\frac{z}{2})^{2m}}{(m!)^2}.$$
Numerically, $C_0(3; 2, 1)\approx 3.66043$ and $C_0(4; 3,1)\approx 3.08214$. When $q$ is large, using the method in \cite{F-M} (Section 2), we can find asymptotic formulas for $C_j(q; a, b)$,
\begin{equation}\nonumber
C_j(q; a, b)=\frac{(2j-1)!! \sqrt{2\pi}}{V(q; a, b)^{j+\frac{1}{2}}}+O_j\(\frac{1}{V(q; a, b)^{j+\frac{3}{2}}} \),
\end{equation}
where $(2j-1)!!=(2j-1)(2j-3)\cdots 3\cdot 1$, $(-1)!!=1$, and
\begin{equation}\nonumber
V(q; a, b)=\sum_{\chi \bmod q} |\chi(b)-\chi(a)|^2 \sum_{\substack{\gamma_{\chi}\in \mathbb{R}\\ L(\frac{1}{2}+i\gamma_{\chi}, \chi)=0}}\frac{1}{\frac{1}{4}+\gamma^2_{\chi}}.
\end{equation}
By Proposition 3.6 in \cite{F-M}, under ${\rm ERH_q}$, $V(q; a, b)\sim 2\phi(q)\log q$.


\section{Formulas for the associated Dirichlet series and orgin of the bias}\label{formu-F-O-o-mega}

Let $\chi$ be a non-principal Dirichlet character, and denote
\begin{equation*}
	F_{f_k}(s,\chi):=\sum_{f(n)=k} \frac{\chi(n)}{n^s},
\end{equation*} 
where $f=\Omega$ or $\omega$. The formulas for $F_{f_k}(s, \chi)$ are needed to analyze the character sums in \eqref{Delta-Omega} and \eqref{Delta-omega}. The purpose of this section is to express $F_{f_k}(s, \chi)$ in terms of Dirichlet $L$-functions, and to explain the source of the biases in the functions $\Delta_{\Omega_k}(x; q, a, b)$ and $\Delta_{\omega_k}(x; q, a, b)$. 

Throughout the paper, the notation $\log z$ will always denote the \textbf{principal branch} of the logarithm of a complex number $z$. 

\subsection{Symmetric functions}

Let $x_1$, $x_2$, $\dots$ be an infinite collection of indeterminates. We say a formal power series  $P(x_1, x_2, \dots)$ with bounded degree is a \textit{symmetric function} if it is invariant under all finite permutations of the variables $x_1$, $x_2$, $\dots$. 

The $n$-th \textit{elementary symmetric function} $e_n=e_n(x_1, x_2, \dots)$ is defined by the generating function
$\sum_{n=0}^{\infty} e_n z^n=\prod_{i=1}^{\infty} (1+x_i z).$
Thus, $e_n$ is the sum of all square-free monomials of degree $n$. Similary, the $n$-th \textit{homogeneous symmetric function} $h_n=h_n(x_1, x_2, \dots)$ is defined by the generating function
$\sum_{n=0}^{\infty} h_n z^n=\prod_{i=1}^{\infty}\frac{1}{1-x_i z}.$
We see that, $h_n$ is the sum of all possible monomials of degree $n$. And the $n$-th \textit{power symmetric function} $p_n=p_n(x_1, x_2, \dots)$ is defined to be
$p_n=x_1^n+x_2^n+\cdots.$

The following result is due to Newton or Girard (see \cite{Mac}, Chapter 1, (2.11) and (2.11'), page 23, or \cite{Me-Re}, Chapter 2, Theorems 2.8 and 2.9).
\begin{lem}\label{lem-Newton}
	For any integer $k\geq 1$,
	\begin{equation}
	kh_k=\sum_{n=1}^k h_{k-n}p_n,\label{lem-Newton-1}
	\end{equation}
	\begin{equation}
	ke_k=\sum_{n=1}^k (-1)^{n-1} e_{k-n}p_n.\label{lem-Newton-2}
	\end{equation}
\end{lem}

\subsection{Formula for $F_{\Omega_k}(s, \chi)$}
For $\Re(s)>1$, we define
\begin{equation}\nonumber
	F(s,\chi):=\sum_{p} \frac{\chi(p)}{p^s},
\end{equation}
the sum being over all prime $p$. Since
\begin{equation}\label{log-L-expr}
\log L(s, \chi)=\sum_{m=1}^{\infty}\sum_{p~\text{prime}} \frac{\chi(p^m)}{m p^{ms}},
\end{equation}
we then have
\begin{equation}\label{F-to-L}
	F(s,\chi)=\log L(s, \chi)-\frac{1}{2}\log L(2s, \chi^2)+G(s),
\end{equation}
where $G(s)$ is absolutely convergent for $\Re(s)\geq \sigma_0$ for any fixed $\sigma_0>\frac{1}{3}$. Henceforth, $\sigma_0$ will be a fixed abscissa $>\frac{1}{3}$, say $\sigma_0=0.34$.  Because $L(s, \chi)$ is an entire function for non-principal characters $\chi$, formula \eqref{F-to-L} provides an analytic continuation of $F(s, \chi)$ to any simply-connected domain within the half-plane $\{s: \Re(s)\geq \sigma_0  \}$ which avoids the zeros of $L(s, \chi)$ and the zeros and possible pole of $L(2s, \chi^2)$. 

For any complex number $s$ with $\Re(s)\geq \sigma_0>\frac{1}{3}$, let $x_p=\frac{\chi(p)}{p^s}$ if $p$ is a prime, $0$ otherwise. Then, by (\ref{lem-Newton-1}) in Lemma \ref{lem-Newton}, we have the following relation
\begin{equation}\label{F-Omega-formu-in}
	kF_{\Omega_k}(s, \chi)=\sum_{n=1}^k F_{\Omega_{k-n}}(s, \chi)F(ns, \chi^n).
\end{equation}
For example, for $k=1$, $F_{\Omega_1}(s, \chi)=F(s, \chi)$. 	For $k=2$,
$$2F_{\Omega_2}(s, \chi)=F^2(s, \chi)+F(2s, \chi^2).$$
For $k=3$,
\begin{align*}
	3! F_{\Omega_3}(s, \chi)&= 2F_{\Omega_2}(s,\chi) F(s,\chi)+2F(s,\chi)F(2s, \chi^2)+2F(3s, \chi^3)\nonumber\\
	&= F^3(s, \chi)+3F(s,\chi)F(2s, \chi^2)+2F(3s, \chi^3).
\end{align*}
For $k=4$,
\begin{align*}
	4! F_{\Omega_4}(s, \chi)&= 3! F_{\Omega_3}(s, \chi)F(s, \chi)+3!F_{\Omega_2}(s,\chi)F(2s, \chi^2)+3!F(s, \chi)F(3s, \chi^2)+3!F(4s, \chi^4)\nonumber\\
	&=F^4(s, \chi)+6F^2(s, \chi)F(2s, \chi^2)+8F(s,\chi)F(3s, \chi^3)+6F(4s, \chi^4)\nonumber\\
	&\quad +3F^2(2s, \chi^2).
\end{align*}

For any integer $l\geq 1$, we define the set
\begin{equation*}
	S_{m,l}^{(k)}:=\{(n_1, \cdots, n_l) ~|~ n_1+\cdots+n_l=k-m, 2\leq n_1\leq n_2\leq \cdots\leq n_l, n_j\in\mathbb{N} (1\leq j\leq l) \}
\end{equation*}
Let $S_m^{(k)}=\bigcup_{l\geq 1} S_{m, l}^{(k)}.$
Thus any element of $S_m^{(k)}$ is a partition of $k-m$ with each part $\geq 2$. For any $\boldsymbol{n}=(n_1, n_2, \cdots, n_l)\in S_m^{(k)}$, denote
\begin{equation*}
	F(\boldsymbol{n}s, \chi):=\prod_{j=1}^l F(n_j s, \chi^{n_j}).
\end{equation*}
Hence, by (\ref{F-Omega-formu-in}) and induction on $k$, we deduce the following result. 
\begin{lem}\label{lem-F-Omega-k}
For $k=1$, $F_{\Omega_1}(s, \chi)=F(s, \chi)$. 	
For any $k\geq 2$, we have
\begin{equation}\label{F-Omega-k}
	k! F_{\Omega_k}(s, \chi)=F^k(s, \chi)+\sum_{m=0}^{k-2}F^m(s, \chi) F_{\boldsymbol{n_m}}(s, \chi),
\end{equation}
where
	$F_{\boldsymbol{n_m}}(s, \chi)=\sum\limits_{\boldsymbol{n}\in S_m^{(k)}} a^{(k)}_m (\boldsymbol{n}) F(\boldsymbol{n}s, \chi)$
for some $a^{(k)}_m (\boldsymbol{n})\in \mathbb{N}$.
\end{lem}

\subsection{Formula for $F_{\omega_k}(s, \chi)$}
By definition, we have
\begin{equation*}
F_{\omega_{k}}(s, \chi)=\sum_{\substack{p_1<p_2<\cdots <p_k\\ p_i~\text{prime}}} \prod_{n=1}^k \(\sum_{j=1}^{\infty} \frac{\chi(p_n^j)}{p_n^j}\).
\end{equation*}
Denote
\begin{equation*}
\widetilde{F}(s, \chi):=\sum_{p~\text{prime}} \(\frac{\chi(p)}{p^s}+\frac{\chi(p^2)}{p^{2s}}+\cdots \),
\end{equation*}
and for any $u\in \mathbb{N}^{+}$,
\begin{equation*}
\widetilde{F}(s, \chi; u):=\sum_{p~\text{prime}} \(\frac{\chi(p)}{p^s}+\frac{\chi(p^2)}{p^{2s}}+\cdots \)^u=\sum_{p~\text{prime}} \sum_{j=u}^{\infty} \(D_u(j)\frac{\chi(p^j)}{p^{js}} \),
\end{equation*}
where $D_u(j)={j-1 \choose u-1}$ is the number of ways of writing $j$ as sum of $u$ ordered positive integers.

By (\ref{log-L-expr}), we have
\begin{equation}\label{F-omega-L}
\widetilde{F}(s, \chi)=\widetilde{F}(s, \chi; 1)=\sum_{p~\text{prime}}\sum_{j=1}^{\infty} \frac{\chi(p^j)}{p^{js}}=\log L(s, \chi)+\frac{1}{2}\log L(2s,\chi^2)+\widetilde{G}_1(s),
\end{equation}
and
\begin{equation}\label{F-omega-L-2s}
\widetilde{F}(s, \chi; 2)=\sum_{p~\text{prime}}\sum_{j=2}^{\infty} (j-1) \frac{\chi(p^j)}{p^{js}}=\log L(2s, \chi^2)+\widetilde{G}_2(s),
\end{equation}
where $\widetilde{G}_1(s)$ and $\widetilde{G}_2(s)$ are absolutely convergent for $\Re(s)\geq \sigma_0$. Formula \eqref{F-omega-L} provides an analytic continuation of $\widetilde{F}(s, \chi)$ to any simply-connected domain within the half-plane $\{s: \Re(s)\geq \sigma_0  \}$ which avoids the zeros of $L(s, \chi)$ and the zeros and possible pole of $L(2s, \chi^2)$. 
Moreover, for any fixed $u\geq 3$, $\widetilde{F}(s, \chi; u)$ is absolutely convergent for $\Re(s)\geq \sigma_0$.

For any complex number $s$ with $\Re(s)\geq \sigma_0$, take $x_p=\sum_{j=1}^{\infty}\frac{\chi(p^j)}{p^{js}}$ if $p$ is a prime, $0$ otherwise. Then by (\ref{lem-Newton-2}) in Lemma \ref{lem-Newton}, we get the following formula,
\begin{equation}\label{F-omega-formu-in}
kF_{\omega_k}(s, \chi)=F_{\omega_{k-1}}(s, \chi)\widetilde{F}(s, \chi)-\sum_{n=2}^k (-1)^n F_{\omega_{k-n}}(s, \chi)\widetilde{F}(s, \chi; n).
\end{equation}
For example, for $k=1$, $F_{\omega_1}(s, \chi)=\widetilde{F}(s, \chi)$. 	For  $k=2$,
\begin{equation*}
2F_{\omega_2}(s, \chi)=\widetilde{F}^2(s, \chi)-\widetilde{F}(s, \chi; 2).
\end{equation*}
For $k=3$,
\begin{align*}
3! F_{\omega_3}(s, \chi)&=2F_{\omega_2}(s,\chi)\widetilde{F}(s, \chi)-2F_{\omega_1}(s, \chi)\widetilde{F}(s, \chi; 2)+2\widetilde{F}(s, \chi; 3)\nonumber\\
&=\widetilde{F}^3(s, \chi)-3\widetilde{F}(s, \chi)\widetilde{F}(s, \chi; 2)+2\widetilde{F}(s, \chi; 3).
\end{align*}
For $k=4$,
\begin{align*}
4! F_{\omega_4}(s, \chi)&=3! F_{\omega_3}(s, \chi)\widetilde{F}(s, \chi)-3!F_{\omega_2}(s, \chi)\widetilde{F}(s, \chi; 2)+3!\widetilde{F}(s, \chi)\widetilde{F}(s, \chi; 3)-3!\widetilde{F}(s, \chi; 4)\nonumber\\
&= \widetilde{F}^4(s, \chi)-6\widetilde{F}^2(s, \chi)\widetilde{F}(s, \chi; 2)+8\widetilde{F}(s, \chi)\widetilde{F}(s, \chi; 3)-6\widetilde{F}(s, \chi; 4)+3\widetilde{F}^2(s, \chi; 2).
\end{align*}

Hence, by (\ref{F-omega-formu-in}) and induction on $k$, we get the following result. 
\begin{lem}\label{lem-F-omega-k}
For $k=1$, $F_{\omega_1}(s, \chi)=\widetilde{F}(s, \chi)$. 	
For any $k\geq 2$, we have
\begin{equation}\label{F-omega-k}
k!F_{\omega_{k}}(s, \chi)= \widetilde{F}^k(s, \chi)+\sum_{m=0}^{k-2} \widetilde{F}^m(s, \chi) \widetilde{F}_{\boldsymbol{n_m}}(s, \chi),
\end{equation}
where
$\widetilde{F}_{\boldsymbol{n_m}}(s, \chi)=\sum\limits_{\boldsymbol{n}\in S_m^{(k)}} b_m^{(k)} (\boldsymbol{n}) \widetilde{F}(\boldsymbol{n}s,\chi)$
for some $b_m^{(k)}(\boldsymbol{n})\in \mathbb{Z}$, and for any $\boldsymbol{n}=(n_1, \cdots, n_l)\in S^{(k)}_m$,
$\widetilde{F}(\boldsymbol{n}s, \chi):=\prod_{j=1}^l \widetilde{F}(s, \chi; n_j).$
\end{lem}


\subsection{Origin of the bias}

In this section, we heuristically explain the origin of the bias in our theorems. 

\textbf{1) Analytical aspect.} In order to get formulas for $\Delta_{\Omega_k}(x; q, a, b)$ and $\Delta_{\omega_k}(x; q, a, b)$, our strategy is to apply Perron's formula to the associated Dirichlet series $F_{\Omega_k}(s, \chi)$ and $F_{\omega_k}(s, \chi)$, then we choose special contours to avoid the singularities of these Dirichlet series. See Section \ref{sec-lems} for the details. 

First, we have a look at the case of counting primes in arithmetic progressions. If we only count primes, by \eqref{F-to-L}, we have
$$F_{\Omega_1}(s, \chi)=F(s, \chi)=\sum_{p}\frac{\chi(p)}{p^s}=\log L(s, \chi)-\frac{1}{2}\log L(2s, \chi^2)+G(s).$$
The main contributions for $\Delta_{\Omega_1}(x; q, a, b)$ are from the first two terms, 
$$\log L(s, \chi)-\frac{1}{2}\log L(2s, \chi^2).$$
The first term $\log L(s, \chi)$ counts all the primes with weight 1 and prime squares with weight $\frac{1}{2}$. The higher order powers of primes are negligible since they only contribute $O(x^{\frac{1}{3}})$.   The singularities of $\log L(s, \chi)$, i.e. the zeros of $L(s, \chi)$, on the critical line contribute the oscillating terms in our result. In our proof, we use special Hankel contours to avoid the singularities of $\log L(s, \chi)$ and extract these oscillating terms (Lemma \ref{lem-main-term}). See Sections \ref{sec-lems} and \ref{Sec-pf} for the details of how to handle these singularities.  
The second term $-\frac{1}{2}\log L(2s, \chi^2)$  counts the prime squares with weight $-\frac{1}{2}$ and contributes the bias term.  When $\chi$ is a real character,  the point $s=\frac{1}{2}$ is a pole of $L(2s, \chi^2)$, and hence the integration of $ -\frac{1}{2}\log L(2s, \chi^2)$ over the Hankel contour around $s=\frac{1}{2}$ contributes a bias term with order of magnitude $\frac{\sqrt{x}}{\log x}$.   Using the orthogonality of  Dirichlet characters, and the formula $\sum_{\chi ~\rm real}(\overline{\chi}(a)-\overline{\chi}(b))=N(q, a)-N(q, b)$, we get the expected size of the bias.  

Another natural and convenient function to consider is $-\frac{L'(s, \chi)}{L(s, \chi)}=\sum_{n=1}^{\infty}\frac{\chi(n)\Lambda(n)}{n^s}$, which is much easier to analyze than $\log L(s, \chi)$. This weighted form is counting each prime $p$ and its powers with weight $\log p$. 
 Similar to $\log L(s, \chi)$, all the singularities of the function $-\frac{L'(s, \chi)}{L(s, \chi)}$ on the critical line are the non-trivial zeros of $L(s, \chi)$ and thus there is no bias for this weighted counting function 
$$\sum_{\substack{n\leq x\\ n\equiv a \bmod q}} \Lambda(n)-\sum_{\substack{n\leq x\\ n\equiv b \bmod q}} \Lambda(n).$$
Thus, partial summation is used to extract the sum 
$$\sum_{\substack{n\leq x\\ n\equiv a \bmod q}}\frac{\Lambda(n)}{\log n}-\sum_{\substack{n\leq x\\ n\equiv b \bmod q}}\frac{\Lambda(n)}{\log n}$$ from the above weighted form, which is possible because $\log n$ is a smooth function.  However, there is no way to do this with the analogue \eqref{Lambda-analog} to recover the unweighted counting function $\Delta_{\Omega_k}(x; q, a, b)$ or $\Delta_{\omega_k}(x; q, a, b)$.

If we count all the prime powers with the same weight 1, by \eqref{F-omega-L}, we have
$$ F_{\omega_1}(s, \chi)=\widetilde{F}(s, \chi)=\log L(s, \chi)+\frac{1}{2}\log L(2s,\chi^2)+\widetilde{G}_1(s). $$
In this case,  the bias is from the second term $\frac{1}{2}\log L(2s,\chi^2)$ for real character $\chi$ which counts the prime squares with positive weight $\frac{1}{2}$. This is why the bias is opposite to the case of counting only primes.

For the general case, when we derive the formula for $\Delta_{\Omega_k}(x; q, a, b)$ using analytic methods, by \eqref{F-Omega-k} in Lemma \ref{lem-F-Omega-k},  the main contributions for $F_{\Omega_k}(s, \chi)$ will be from $\frac{1}{k!} F^k(s, \chi)$, which is essentially 
$$\frac{1}{k!}\(\log L(s, \chi)-\frac{1}{2}\log L(2s, \chi^2)\)^k.$$
In the expansion of the above formula, the term $\frac{1}{k!}\log^k L(s, \chi)$ contributes the oscillating terms (see \eqref{rho-form} and \eqref{main-rho})
$$\frac{(-1)^k}{(k-1)!}\frac{\sqrt{x}(\log\log x)^{k-1}}{\log x}\sum_{L(\frac{1}{2}+i\gamma_{\chi}, \chi)=0} \frac{x^{i\gamma_{\chi}}}{\frac{1}{2}+i\gamma_{\chi}}.$$
When $\chi$ is real, the term 
$$\frac{1}{k!}\(-\frac{1}{2}\log L(2s, \chi^2)\)^k=\frac{(-1)^k}{k! 2^k}\(\log L(2s, \chi^2)\)^k$$ 
contributes a bias term (see \eqref{bias-form} and \eqref{main-bias})
$$\frac{1}{(k-1)!}\frac{(-1)^k}{ 2^{k-1}}\frac{\sqrt{x}(\log\log x)^{k-1}}{\log x}.$$
Then summing over all the real characters, we get the expected bias term in our formula for $\Delta_{\Omega_k}(x; q, a, b)$. The factor $\frac{(-1)^k}{2^{k-1}}$ explains why the bias has different directions depending on the parity of $k$ and why the bias decreases as $k$ increases. 
Other terms with factors of the form $\log^{k-j} L(s, \chi) \log^j (2s, \chi^2)$ for $1\leq j\leq k-1$ only contribute oscillating terms with lower orders of $\log\log x$ which can be put into the error term in our formula (see Lemma \ref{lem-distri-sum}). 

Similarly, for the case of $\Delta_{\omega_k}(x; q, a, b)$, by \eqref{F-omega-k} in Lemma \ref{lem-F-omega-k}, the main contributions for $F_{\omega_k}(s, \chi)$ are from $$\frac{1}{k!}\widetilde{F}^k(s, \chi)=\frac{1}{k!}\(\log L(s, \chi)+\frac{1}{2}\log L(2s,\chi^2)+\widetilde{G}_1(s) \)^k.$$
The main terms are from the contributions of the terms $\frac{1}{k!}\log^k L(s, \chi)$ and $\frac{1}{k!} \(\frac{1}{2}\log L(2s, \chi^2)\)^k$. Thus, the main oscillating terms are the same as that of $\Delta_{\Omega_k}(x; q, a, b)$, and the bias term has the same size without direction change. 

Through the above analysis, we see that the biases are mainly affected by the powers of $\pm\frac{1}{2}\log L(2s, \chi^2)$ for real characters which count the products of prime squares. 

\bigskip

 \textbf{2) Combinatorial aspect.}  Instead of giving precise prediction of the size of the bias as above, here we use a simpler combinatorial intuition to roughly explain the behavior of the bias. We borrowed this combinatorial explanation from Hudson \cite{Hudson}. 

Pick a large number $X$. Let $S_1$ be the set of primes $p\equiv 1 \bmod 4$ up to $X$, and $S_2$ be the set of primes $p\equiv 3\bmod 4$ up to $X$. Using these primes, we  generate the set $V^{(2)}:=\{ pq: p,  q \in S_1 \cup S_2, p \text{~and~} q \text{~can be the same}\}$.  

Let $V_1^{(2)}:=\{n\in V^{(2)}: n\equiv 1\bmod 4 \}$, and $V_2^{(2)}:=\{n\in V^{(2)}: n\equiv 3\bmod 4 \}$. Then, the integers in $V_1^{(2)}$ come from either products of two primes from $S_1$ or products of two primes from $S_2$.  The integers in $V^{(2)}_2$ are the product of two primes $pq$ with $p\in S_1$ and $q\in S_2$. Thus, 
$$|V^{(2)}_1|={|S_1|\choose 2}+|S_1|+{|S_2|\choose 2}+|S_2|=\frac{|S_1|^2+|S_2|^2|}{2}+\frac{|S_1|+|S_2|}{2},$$
and 
$$|V^{(2)}_2|=|S_1|\cdot |S_2|.$$
It is clear that $|V^{(2)}_1|>|V^{(2)}_2|$. Note that $\frac{|S_1|+|S_2|}{2}$ counts the squares of primes with weight $\frac{1}{2}$ which makes a crucial difference between $V_1^{(2)}$ and $V_2^{(2)}$. 

Let $V_1^{(0)}=\{ 1\}$ and $V_2^{(0)}=\emptyset$. For any $k\geq 1$, denote 
$$V_1^{(k)}:=\{n=p_1\cdots p_k:  p_j\in S_1\cup S_2 \text{~for all~} 1\leq j\leq k,  n\equiv 1\bmod 4\},$$
$$V_2^{(k)}:=\{n=p_1\cdots p_k:  p_j\in S_1\cup S_2 ~\text{for all~} 1\leq j\leq k,  n\equiv 3\bmod 4\},$$
where the $p_j$ can be the same. Note that $V_1^{(1)}=S_1$ and $V_2^{(1)}=S_2$. 

We give inductive formulas for $|V_1^{(k)}|$ and $|V_2^{(k)}|$.  The elements of $V_1^{(k)}$ and $V_2^{(k)}$ are generated by integers of the form $ p^j n_{k-j}$ for $p\in S_1 \text{~or~} S_2$ and $n_{k-j}\in V_1^{(k-j)} \text{~or~} V_2^{(k-j)} ~(1\leq j\leq k)$.  By \eqref{lem-Newton-1} in Lemma \ref{lem-Newton}, we have
\begin{align*}
k|V_1^{(k)}|=&\underbrace{\(|V_1^{(k-1)}|\cdot |S_1|+|V_2^{(k-1)}|\cdot |S_2|\)}_{p n_{k-1}}+  \underbrace{|V_1^{(k-2)}| (|S_1|+|S_2|)}_{p^2 n_{k-2}}\\
&+\underbrace{\(|V_1^{(k-3)}|\cdot |S_1|+|V_2^{(k-3)}|\cdot |S_2|\)}_{p^3 n_{k-3}}+\cdots 
\end{align*} 
and 
\begin{align*}
k|V_2^{(k)}|=&\underbrace{\(|V_2^{(k-1)}|\cdot |S_1|+|V_1^{(k-1)}|\cdot |S_2|\)}_{p n_{k-1}}+  \underbrace{|V_2^{(k-2)}| (|S_1|+|S_2|)}_{p^2 n_{k-2}}\\
&+\underbrace{\(|V_2^{(k-3)}|\cdot |S_1|+|V_1^{(k-3)}|\cdot |S_2|\)}_{p^3 n_{k-3}}+\cdots.
\end{align*} 
Thus, 
\begin{align}\label{combi-explain-Vk}
k\( |V_1^{(k)}|-|V_2^{(k)}|\)=&\(|V_1^{(k-1)}|\cdot |S_1|+|V_2^{(k-1)}|\cdot |S_2|\)-\(|V_2^{(k-1)}|\cdot |S_1|+|V_1^{(k-1)}|\cdot |S_2|\)\nonumber\\
&+\( |V_1^{(k-2)}|-|V_2^{(k-2)}|\) (|S_1|+|S_2|)+\(|V_1^{(k-3)}|\cdot |S_1|+|V_2^{(k-3)}|\cdot |S_2|\)\nonumber\\
&-\(|V_2^{(k-3)}|\cdot |S_1|+|V_1^{(k-3)}|\cdot |S_2|\)+\cdots\nonumber\\
=&\( |V_1^{(k-1)}|-|V_2^{(k-1)}|\) ( |S_1|-|S_2|)+\( |V_1^{(k-2)}|-|V_2^{(k-2)}|\) (|S_1|+|S_2|)\nonumber\\
&+\( |V_1^{(k-3)}|-|V_2^{(k-3)}|\) ( |S_1|-|S_2|)+\cdots.
\end{align}

For $1\leq j\leq k-1$, suppose $|V_1^{(j)}|<|V_2^{(j)}|$ for odd $j$ and $|V_1^{(j)}|>|V_2^{(j)}|$ for even $j$. Therefore, by \eqref{combi-explain-Vk} and induction, we deduce that $|V_1^{(k)}|<|V_2^{(k)}|$ for odd $k$ and $|V_1^{(k)}|>|V_2^{(k)}|$ for even $k$. 
This provides us a heuristic explanation for the bias oscillation of $\Delta_{\Omega_k}(x; q, a, b)$.

\section{Contour integral representation}\label{sec-lems}

In this section, we express the inner sums in \eqref{Delta-Omega} and \eqref{Delta-omega} as integrals over  truncated Hankel contours (see Lemma \ref{lem-A} below).

Let
\begin{equation*}
\psi_{f_k}(x,\chi):=\sum_{\substack{n\leq x\\ f(n)=k}}\chi(n),
\end{equation*}
where $f=\Omega$ or $\omega$. 
By Perron's formula  (\cite{Kara}, Chapter V, Theorem 1), we have the following lemma.
\begin{lem}\label{lem-psi-Perron}
	For any $T\geq 2$,
	\begin{equation*}
	\psi_{f_k}(x,\chi)=\frac{1}{2\pi i}\int_{c-iT}^{c+iT} F_{f_k}(s, \chi)\frac{x^s}{s} ds +O\(\frac{x\log x}{T}+1\),
	\end{equation*}
	where $c=1+\frac{1}{\log x}$, and $f=\Omega$ or $\omega$.
\end{lem}

Starting from Lemma \ref{lem-psi-Perron}, we will shift the contour to the left, in a way which avoids the singularities of the integrand. We will then require estimates of the integrand along the various parts of the new contour. 

\begin{lem}\label{lem-F-est}
	Assume $ERH_q$. Then, for any $0<\delta<\frac{1}{6}$ and for all $\chi\neq \chi_0 \bmod q$, there exists a sequence of numbers $\mathcal{T}=\{T_n \}_{n=0}^{\infty}$ satisfying $n\leq T_n\leq n+1$ such that, for $T\in \mathcal{T}$,
	$$F_{f_k}(\sigma+iT)=O\(\log^k T \), ~( \frac{1}{2}-\delta<\sigma<2)$$
	where $f=\Omega$ or $\omega$. 	
\end{lem}
\noindent \textit{\textbf{Proof.}}
Using the similar method as in \cite{Titch} (Theorem 14.16), one can show that, for any $\epsilon>0$ and for all $\chi\neq \chi_0 \bmod q$, there exists a sequence of numbers $\mathcal{T}=\{T_n \}_{n=0}^{\infty}$ satisfying $n\leq T_n\leq n+1$ such that,
$T_n^{-\epsilon}\ll |L(\sigma+iT_n, \chi)|\ll T_n^{\delta+\epsilon}, ~(\frac{1}{2}-\delta<\sigma<2).$
Hence, by formulas (\ref{F-to-L}), (\ref{F-Omega-k}), (\ref{F-omega-L}), (\ref{F-omega-L-2s}), and (\ref{F-omega-k}), we get the conclusion of this lemma. \qed

\bigskip

Let $\rho$ be a zero of $L(s, \chi)$, $\Delta_{\rho}$ be the distance of $\rho$ to the nearest other zero, and $D_{\gamma}:=\min\limits_{T\in \mathcal{T}}(|\gamma-T|)$.
For each zero $\rho$, and $X>0$, let $\mathcal{H}(\rho, X)$ denote the truncated Hankel contour surrounding the point $s=\rho$ with radius $0<r_{\rho}\leq \min(\frac{1}{x}, \frac{\Delta_{\rho}}{3}, \frac{D_{\gamma}}{2}, \frac{|\rho-1/2|}{3})$, which includes the circle $|s-\rho|=r_{\rho}$ excluding the point $s=\rho-r_\rho$, and the half-line $(\rho-X, \rho-r]$ traced twice with arguments $+\pi$ and $-\pi$ respectively. Let $\Delta_0$ be the distance of $\frac{1}{2}$ to the nearest zero. Let $\mathcal{H}(\frac{1}{2}, X)$ denote the corresponding truncated Hankel contour surrounding $s=\frac{1}{2}$ with radius $r_0=\min(\frac{1}{x}, \frac{\Delta_0}{3})$.

\begin{figure}[h!]
	\centering
	\includegraphics[width=10cm]{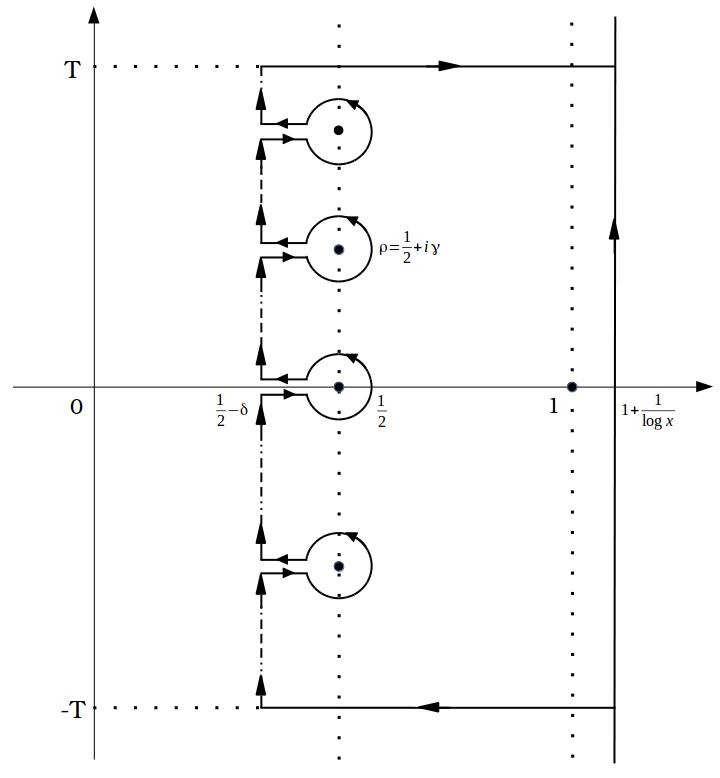}\\
	\caption{Integration contour}\label{contour}
\end{figure}

Take $\delta=\frac{1}{10}$. By Lemma \ref{lem-psi-Perron}, we pull the contour to the left to the line $\Re(s)=\frac{1}{2}-\delta$ using the truncated Hankel contour  $\mathcal{H}(\rho, \delta)$ to avoid the zeros of $L(s, \chi)$ and using $\mathcal{H}(\frac{1}{2}, \delta)$ to avoid the point $s=\frac{1}{2}$. See Figure \ref{contour}.

Then we have the following lemma.

\begin{lem}\label{lem-A}
	Assume $ERH_q$, and $L(\frac{1}{2}, \chi)\neq 0$ $(\chi\neq \chi_0)$. Then, for any fixed $k\geq 1$, and $T\in \mathcal{T}$,
	\begin{align*}
	\psi_{f_k}(x,\chi)&=\sum_{|\gamma|\leq T}\frac{1}{2\pi i}\int_{\mathcal{H}(\rho, \delta)}F_{f_k}(s, \chi)\frac{x^s}{s} ds +a(\chi)\frac{1}{2\pi i}\int_{\mathcal{H}(\frac{1}{2},  {\delta})}F_{f_k}(s, \chi)\frac{x^s}{s} ds\\
	&\quad +O\(\frac{x\log x}{T}+\frac{x(\log T)^k}{T}+x^{\frac{1}{2}-\delta} (\log T)^{k+1}\),
	\end{align*}
	where $a(\chi)=1$ if $\chi$ is real, $0$ otherwise, and $f=\omega$ or $\Omega$.
\end{lem}

\noindent \textbf{\textit{Proof.}}  By formulas (\ref{F-Omega-k}) and (\ref{F-omega-k}), if $\chi$ is not real, $s=\frac{1}{2}$ is not a singularity of $F_{f_k}(s, \chi)$. Hence the second term is zero if $\chi$ is not real. By Lemma \ref{lem-F-est}, the integral on the horizontal line is
 \begin{eqnarray}\label{lem-A-E-1}
 \ll (\log T)^k\int_{\frac{1}{2}-\delta}^{c}\frac{x^{\sigma}}{|\sigma+iT|} d\sigma \ll   \frac{x^c(\log T)^k}{T} \ll \frac{x(\log T)^k}{T}.
 \end{eqnarray}
Under the assumption ${\rm ERH_q}$, the integral on the vertical line $\Re(s)=\frac{1}{2}-\delta$ is
\begin{eqnarray}\label{lem-A-E-2}
\ll \int_{-T}^T \frac{x^{\frac{1}{2}-\delta} \log^k (|t|+2)}{|\frac{1}{2}-\delta+it|}dt\ll x^{\frac{1}{2}-\delta} (\log T)^{k+1}.
\end{eqnarray}
By (\ref{lem-A-E-1}), (\ref{lem-A-E-2}), and Lemma \ref{lem-psi-Perron}, we get the desired error term in this lemma. \qed

\bigskip


\section{Proof of  the main theorems }\label{Sec-pf}

In this section, we outline the proof of Theorems \ref{thm-Delta-Omega} and \ref{thm-Delta-omega}.

Let $\gamma$ be the imaginary part of a zero of $L(s, \chi)$ in the critical strip. We have the following lemma. 
\begin{lem}[\cite{Ford-S}, Lemma 2.2] \label{lem-Ford-S-2.2}
	Let $\chi$ be a Dirichlet character modulo $q$. Let $N(T, \chi)$ denote the number of zeros of $L(s, \chi)$ with $0<\Re(s)<1$ and $|\Im(s)|<T$. Then
	
	1) $N(T, \chi)=O(T\log (qT))$ for $T\geq 1$.
	
	2) $N(T, \chi)-N(T-1, \chi)=O(\log (qT))$ for $T\geq 1$.
	
	3) Uniformly for $s=\sigma+it$ and $\sigma\geq -1$,
	\begin{equation}\label{F-S-lem-formu}
	\frac{L'(s, \chi)}{L(s, \chi)}=\sum_{|\gamma-t|<1} \frac{1}{s-\rho} +O(\log q(|t|+2)).
	\end{equation}
\end{lem}

For simplicity,  we denote
$$\frac{1}{\Gamma_j(u)}:=\left[\frac{d^j}{dz^j}\(\frac{1}{\Gamma(z)}\)\right]_{z=u}.$$
The following lemma is the starting lemma to give us the bias terms and oscillating terms in our main theorems. This lemma may have independent use, we will give the proof in Section \ref{sec-main-lem}.
\begin{lem}\label{lem-main-term}
	Let $\mathcal{H}(a, \delta)$  be the truncated Hankel contour surrounding a complex number $a \ (\Re(a)>2\delta)$ with radius $0<r\ll \frac{1}{x}$. Then, for any integer $k\geq 1$,
	\begin{align*}
	&\frac{1}{2\pi i}\int_{\mathcal{H}(a, \delta)} \log^k(s-a)\frac{x^s}{s}ds\nonumber\\
	&=\frac{(-1)^k x^a}{a\log x}\left\{  k (\log\log x)^ {k-1}+\sum_{j=2}^k {k \choose j} (\log\log x)^{k-j} \frac{1}{\Gamma_j(0)}\right\}+O_k\(\frac{|x^{a-\delta/3}|}{|a|}\)\nonumber\\
	&\quad +O_k\(\frac{|x^a|}{|a|^2 \log^2 x} (\log\log x)^{k-1}\)+O_k\( \frac{|x^a|}{|a|^2|\Re(a)-\delta|}\frac{(\log\log x)^{k-1}}{(\log x)^3}\)
	.
	\end{align*}
\end{lem}
\noindent\textbf{Remark.} By (\ref{pf-lem-Esigma-I-j-l}) in the proof of Lemma \ref{lem-Esigma}, one can easily show that
\begin{equation}\label{Gamma-j-estimate}
\left|\frac{1}{\Gamma_j(0)}\right|\ll \Gamma(j+1).
\end{equation}

By Lemma \ref{lem-A},  we need to examine the integration over the truncated Hankel contours $\mathcal{H}(\rho, \delta)$ and $\mathcal{H}(\frac{1}{2}, \delta)$. 
By (\ref{F-to-L}) and (\ref{F-omega-L}), and the assumptions of our theorems, on each truncated Hankel contour $\mathcal{H}(\rho, \delta)$, we integrate the formula \eqref{F-S-lem-formu} in Lemma \ref{lem-Ford-S-2.2} to obtain
\begin{equation}\label{pf-thm-F-rho-H}
F(s, \chi)=\log(s-\rho)+H_{\rho}(s),
\end{equation}
\begin{equation}\label{pf-thm-F-omega-rho-H}
\widetilde{F}(s, \chi)=\log(s-\rho)+\widetilde{H}_{\rho}(s),
\end{equation}
where
\begin{align*}
H_{\rho}(s)&=\sum_{0<|\gamma'-\gamma|\leq 1} \log(s-\rho')+O(\log |\gamma|),\\
\widetilde{H}_{\rho}(s)&=\sum_{0<|\gamma'-\gamma|\leq 1} \log(s-\rho')+O(\log |\gamma|).
\end{align*}
 If $\chi$ is real, $s=\frac{1}{2}$ is a pole of $L(2s, \chi^2)$. So, by (\ref{F-to-L}) and (\ref{F-omega-L}), on the truncated Hankel contour $\mathcal{H}(\frac{1}{2}, \delta)$, for a real character $\chi$, we write
\begin{align}
F(s, \chi)&=\frac{1}{2}\log \(s-\frac{1}{2}\)+H_B(s),\label{pf-thm-bias-F-formu}\\
\widetilde{F}(s, \chi)&=-\frac{1}{2}\log \(s-\frac{1}{2}\)+\widetilde{H}_B(s),\label{pf-thm-bias-F-omega-formu}
\end{align}
where $H_B(s)=O(1)$ and $\widetilde{H}_B(s)=O(1)$.

Denote
\begin{align*}
I_{\rho}(x)&:=\frac{1}{2\pi i}\int_{\mathcal{H}(\rho,\delta)}k!F_{\Omega_k}(s, \chi)\frac{x^s}{s} ds, \\
I_B(x)&:= \frac{1}{2\pi i}\int_{\mathcal{H}(\frac{1}{2},\delta)}k!F_{\Omega_k}(s, \chi)\frac{x^s}{s} ds,
\end{align*}
and
\begin{align*}
	\widetilde{I}_{\rho}(x)&:=\frac{1}{2\pi i}\int_{\mathcal{H}(\rho,\delta)}k!F_{\omega_k}(s, \chi)\frac{x^s}{s} ds, \\
	\widetilde{I}_B(x)&:= \frac{1}{2\pi i}\int_{\mathcal{H}(\frac{1}{2},\delta)}k!F_{\omega_k}(s, \chi)\frac{x^s}{s} ds.
\end{align*}

We define a function $T(x)$ as follows: for $T_{n'}\in \mathcal{T}$ satisfying $e^{2^{n+1}}\leq T_{n'}\leq e^{2^{n+1}}+1$, let $T(x)=T_{n'}$ for $e^{2^n}\leq x\leq e^{2^{n+1}}$. In particular, we have
$$x\leq T(x)\leq 2x^2 \quad (x\geq e^2).$$

Thus, by Lemma \ref{lem-A}, for $T=T(x)$,
\begin{align}
\psi_{\Omega_k}(x, \chi)&=\frac{1}{k!}\sum_{|\gamma|\leq T}I_{\rho}(x)+\frac{a(\chi)}{k!}I_B(x)+O\(x^{\frac{1}{2}-\frac{\delta}{2}} \),\label{psi-chi-Omega-form}\\
\psi_{\omega_k}(x, \chi)&=\frac{1}{k!}\sum_{|\gamma|\leq T}\widetilde{I}_{\rho}(x)+\frac{a(\chi)}{k!}\widetilde{I}_B(x)+O\(x^{\frac{1}{2}-\frac{\delta}{2}}\).\label{psi-chi-omega-form}
\end{align}
We will see later that $ \sum_{|\gamma|\leq T}I_{\rho}(x)$ and $\sum_{|\gamma|\leq T}\widetilde{I}_{\rho}(x)$ will contribute the oscillating terms, i.e. the summation over zeros, in our theorems, and $I_B(x)$ and $\widetilde{I}_B(x)$ will contribute the bias terms.

Next, we want to find the main contributions for $I_{\rho}(x)$, $I_B(x)$, $\widetilde{I}_{\rho}(x)$, and $\widetilde{I}_B(x)$.
By (\ref{F-Omega-k}) and (\ref{pf-thm-F-rho-H}), we have
\begin{align}\label{rho-form}
I_{\rho}(x)&=\frac{1}{2\pi i}\int_{\mathcal{H}(\rho, \delta)} (\log(s-\rho))^k \frac{x^s}{s} ds+\frac{1}{2\pi i}\int_{\mathcal{H}(\rho, \delta)} \sum_{j=1}^k {k \choose j} \(\log(s-\rho)\)^{k-j} (H_{\rho}(s))^j\frac{x^s}{s}ds\nonumber\\
&\quad +\frac{1}{2\pi i}\int_{\mathcal{H}(\rho, \delta)} \sum_{m=0}^{k-2} F^m(s ,\chi) F_{\boldsymbol{n_m}}(s, \chi)\frac{x^s}{s}ds\nonumber\\
&=: I_{M_\rho}(x)+E_{M_\rho}(x)+E_{R_\rho}(x),
\end{align}
and by (\ref{F-Omega-k}) and (\ref{pf-thm-bias-F-formu}),
\begin{align}\label{bias-form}
I_B(x)&=\frac{1}{2^k} \frac{1}{2\pi i} \int_{\mathcal{H}(\frac{1}{2},\delta)} \(\log\(s-
\frac{1}{2}\)\)^k \frac{x^s}{s}ds\nonumber\\
&\quad +\frac{1}{2\pi i} \int_{\mathcal{H}(\frac{1}{2},\delta)} \sum_{j=1}^k {k \choose j} \(\frac{1}{2}\log\(s-\frac{1}{2}\)   \)^{k-j} \(H_B(s)\)^j \frac{x^s}{s}ds \nonumber\\
&\quad +\frac{1}{2\pi i} \int_{\mathcal{H}(\frac{1}{2},\delta)} \sum_{m=0}^{k-2} F^m(s ,\chi) F_{\boldsymbol{n_m}}(s, \chi)\frac{x^s}{s}ds\nonumber\\
&=: B_M(x)+E_B(x)+E_R(x).
\end{align}
Here, $I_{M_{\rho}}(x)$ and $B_M(x)$ will make main contributions to $I_{\rho}(x)$ and $I_B(x)$, respectively. Similarly, by (\ref{F-omega-k}) and (\ref{pf-thm-F-omega-rho-H}), we have
\begin{align}\label{rho-form-omega}
\widetilde{I}_{\rho}(x)&=\frac{1}{2\pi i}\int_{\mathcal{H}(\rho, \delta)} (\log(s-\rho))^k \frac{x^s}{s} ds+\frac{1}{2\pi i}\int_{\mathcal{H}(\rho, \delta)} \sum_{j=1}^k {k \choose j} \(\log(s-\rho)\)^{k-j} (\widetilde{H}_{\rho}(s))^j\frac{x^s}{s}ds\nonumber\\
&\quad +\frac{1}{2\pi i}\int_{\mathcal{H}(\rho, \delta)} \sum_{m=0}^{k-2} \widetilde{F}^m(s ,\chi) \widetilde{F}_{\boldsymbol{n_m}}(s, \chi)\frac{x^s}{s}ds\nonumber\\
&=: \widetilde{I}_{M_\rho}(x)+\widetilde{E}_{M_\rho}(x)+\widetilde{E}_{R_\rho}(x),
\end{align}
and by (\ref{F-omega-k}) and (\ref{pf-thm-bias-F-omega-formu}),
\begin{align}\label{bias-form-omega}
\widetilde{I}_B(x)&=\frac{(-1)^k}{2^k} \frac{1}{2\pi i} \int_{\mathcal{H}(\frac{1}{2},\delta)} \(\log\(s-
\frac{1}{2}\)\)^k \frac{x^s}{s}ds\nonumber\\
&\quad +\frac{1}{2\pi i} \int_{\mathcal{H}(\frac{1}{2},\delta)} \sum_{j=1}^k {k \choose j} \(-\frac{1}{2}\log\(s-\frac{1}{2}\)   \)^{k-j} \(\widetilde{H}_B(s)\)^j \frac{x^s}{s}ds \nonumber\\
&\quad +\frac{1}{2\pi i} \int_{\mathcal{H}(\frac{1}{2},\delta)} \sum_{m=0}^{k-2} \widetilde{F}^m(s ,\chi) \widetilde{F}_{\boldsymbol{n_m}}(s, \chi)\frac{x^s}{s}ds\nonumber\\
&=: \widetilde{B}_M(x)+\widetilde{E}_B(x)+\widetilde{E}_R(x).
\end{align}
Here, $\widetilde{I}_{M_{\rho}}(x)$ and $\widetilde{B}_M(x)$ will make main contributions to $\widetilde{I}_{\rho}(x)$ and $\widetilde{I}_B(x)$, respectively.

Applying Lemma \ref{lem-main-term}, we have
\begin{align}\label{main-rho}
I_{M_{\rho}}(x)=\widetilde{I}_{M_{\rho}}(x)
&=\frac{(-1)^k\sqrt{x}}{\log x}\frac{x^{i\gamma}}{\frac{1}{2}+i\gamma}\left\{  k (\log\log x)^ {k-1}+\sum_{j=2}^k {k \choose j} (\log\log x)^{k-j} \frac{1}{\Gamma_j(0)}\right\}\nonumber\\
&\quad +O_k\( \frac{1}{|\gamma|^2} \frac{\sqrt{x}(\log\log x)^{k-1}}{(\log x)^2}\)+O_k\(\frac{x^{\frac{1}{2}-\frac{\delta}{3}}}{|\gamma|}\),
\end{align}

\begin{align}\label{main-bias}
B_M(x)&=\frac{(-1)^k\sqrt{x}}{2^{k-1}\log x}\left\{  k (\log\log x)^ {k-1}+\sum_{j=2}^k {k \choose j} (\log\log x)^{k-j} \frac{1}{\Gamma_j(0)}\right\}\nonumber\\
&\quad +O_k\( \frac{\sqrt{x}(\log\log x)^{k-1}}{(\log x)^2}\)+O_k\(x^{\frac{1}{2}-\frac{\delta}{3}}\),
\end{align}
and
\begin{equation}\label{main-bias-omega}
\widetilde{B}_M(x)=(-1)^k B_M(x).
\end{equation}

For the bias terms, by (\ref{bias-form}), (\ref{bias-form-omega}), (\ref{main-bias}), and (\ref{main-bias-omega}), we have
\begin{align}
I_B(x)&=\frac{(-1)^k\sqrt{x}}{2^{k-1}\log x}\left\{  k (\log\log x)^ {k-1}+\sum_{j=2}^k {k \choose j} (\log\log x)^{k-j} \frac{1}{\Gamma_j(0)}\right\}\nonumber\\
&\quad +E_B(x)+E_R(x)+O_k\( \frac{\sqrt{x}(\log\log x)^{k-1}}{(\log x)^2}\)+O_k\(x^{\frac{1}{2}-\frac{\delta}{3}}\),
\end{align}
and 
\begin{align}
\widetilde{I}_B(x)&=\frac{\sqrt{x}}{2^{k-1}\log x}\left\{  k (\log\log x)^ {k-1}+\sum_{j=2}^k {k \choose j} (\log\log x)^{k-j} \frac{1}{\Gamma_j(0)}\right\}\nonumber\\
&\quad +\widetilde{E}_B(x)+\widetilde{E}_R(x)+O_k\( \frac{\sqrt{x}(\log\log x)^{k-1}}{(\log x)^2}\)+O_k\(x^{\frac{1}{2}-\frac{\delta}{3}}\).
\end{align}
We will prove the following result in Section \ref{subsec-bias}. 
\begin{lem}\label{lem-bias}
	For the bias terms,
	\begin{equation*}
	I_B(x)=\frac{(-1)^k k}{2^{k-1}}\frac{\sqrt{x}}{\log x} (\log\log x)^{k-1}+O_k\( \frac{\sqrt{x}(\log\log x)^{k-2}}{\log x}\),
	\end{equation*}
	\begin{equation*}
	\widetilde{I}_B(x)=\frac{ k}{2^{k-1}}\frac{\sqrt{x}}{\log x} (\log\log x)^{k-1}+O_k\( \frac{\sqrt{x}(\log\log x)^{k-2}}{\log x}\).
	\end{equation*}
\end{lem}

Then for the oscillating terms, by (\ref{rho-form}), (\ref{rho-form-omega}), and (\ref{main-rho}), and Lemma \ref{lem-Ford-S-2.2}, for $T=T(x)$,
\begin{align}\label{sum-rho-form}
\sum_{|\gamma|\leq T} I_{\rho}(x)
&=\frac{(-1)^k k\sqrt{x}(\log\log x)^ {k-1}}{\log x} \sum_{|\gamma|\leq T}\frac{x^{i\gamma}}{\frac{1}{2}+i\gamma}+\frac{(-1)^k \sqrt{x}}{\log x}\sum_{j=2}^k {k \choose j}  \frac{(\log\log x)^{k-j}}{\Gamma_j(0)}\sum_{|\gamma|\leq T}\frac{x^{i\gamma}}{\frac{1}{2}+i\gamma}\nonumber\\
&\quad + \sum_{|\gamma|\leq T} E_{M_{\rho}}(x) + \sum_{|\gamma|\leq T} E_{R_{\rho}}(x)+O_k\(\frac{\sqrt{x}(\log\log x)^{k-1}}{\log^2 x} \),
\end{align}
and
\begin{align}\label{sum-rho-form-omega}
\sum_{|\gamma|\leq T} \widetilde{I}_{\rho}(x)&=\frac{(-1)^k k\sqrt{x}(\log\log x)^ {k-1}}{\log x} \sum_{|\gamma|\leq T}\frac{x^{i\gamma}}{\frac{1}{2}+i\gamma}+\frac{(-1)^k \sqrt{x}}{\log x}\sum_{j=2}^k {k \choose j}  \frac{(\log\log x)^{k-j}}{\Gamma_j(0)}\sum_{|\gamma|\leq T}\frac{x^{i\gamma}}{\frac{1}{2}+i\gamma}\nonumber\\
&\quad + \sum_{|\gamma|\leq T} \widetilde{E}_{M_{\rho}}(x) + \sum_{|\gamma|\leq T} \widetilde{E}_{R_{\rho}}(x)+O_k\(\frac{\sqrt{x}(\log\log x)^{k-1}}{\log^2 x} \).
\end{align}

The first terms in the above formulas  are  the main oscillating terms in our theorems. We will show in Section \ref{sec-error-terms} that the other terms are small in average. 
For $T=T(x)$, denote
\begin{align}
	\Sigma_1(x; \chi)&:=\frac{\log x}{\sqrt{x}} \sum_{|\gamma|\leq T} E_{M_{\rho}}(x)=\log x \sum_{|\gamma|\leq T} x^{i\gamma} E'_{M_{\rho}}(x), \label{de-Sigma-1}\\
	\Sigma_2(x; \chi)&:=\frac{\log x}{\sqrt{x}} \sum_{|\gamma|\leq T} E_{R_{\rho}}(x)= \log x \sum_{|\gamma|\leq T} x^{i\gamma} E'_{R_{\rho}}(x), \label{de-Sigma-2}
\end{align}
where $E'_{M_{\rho}}(x)=\frac{E_{M_{\rho}}(x)}{x^{\rho}}$, and $E'_{R_{\rho}}(x)=\frac{E_{R_{\rho}}(x)}{x^{\rho}}$. Similarly, denote
\begin{align}
	\widetilde{\Sigma}_1(x; \chi)&:=\frac{\log x}{\sqrt{x}} \sum_{|\gamma|\leq T} \widetilde{E}_{M_{\rho}}(x)= \log x \sum_{|\gamma|\leq T} x^{i\gamma} \widetilde{E}'_{M_{\rho}}(x), \label{de-Sigma-1-omega}\\
	\widetilde{\Sigma}_2(x; \chi)&:=\frac{\log x}{\sqrt{x}} \sum_{|\gamma|\leq T} \widetilde{E}_{R_{\rho}}(x)= \log x \sum_{|\gamma|\leq T} x^{i\gamma} \widetilde{E}'_{R_{\rho}}(x), \label{de-Sigma-2-omega}
\end{align}
where $\widetilde{E}'_{M_{\rho}}(x)=\frac{\widetilde{E}_{M_{\rho}}(x)}{x^{\rho}}$, and $\widetilde{E}'_{R_{\rho}}(x)=\frac{\widetilde{E}_{R_{\rho}}(x)}{x^{\rho}}$.

\bigskip

Then we have the following lemma (see  Section \ref{subsec-error-Hankel} and Section \ref{subsec-esti-Hankel} for the proof). 
\begin{lem}\label{lem-distri-sum} For the error terms from the Hankel contours around zeros, we have
	\begin{align*}
	\int_2^Y \(\left| \Sigma_1(e^y; \chi)\right|^2+\left| \Sigma_2(e^y; \chi)\right|^2 \) dy&= o\(Y(\log Y)^{2k-2}\),\\
	\int_2^Y \(\left| \widetilde{\Sigma}_1(e^y; \chi)\right|^2+\left| \widetilde{\Sigma}_2(e^y; \chi)\right|^2 \) dy&= o\(Y(\log Y)^{2k-2}\).
	\end{align*}
\end{lem}

Moreover, we also need to bound the lower order sum 
\begin{equation}\label{pf-distri-S-1}
S_1(x; \chi):=(-1)^k \sum_{j=2}^k {k \choose j} (\log\log x)^{k-j} \frac{1}{\Gamma_j(0)}\sum_{|\gamma|\leq T}\frac{x^{i\gamma}}{\frac{1}{2}+i\gamma},
\end{equation}
and the error  from the truncation by a fixed large $T_0$,
\begin{equation}\label{pf-distri-S-2}
S_2(x, T_0; \chi):=\sum_{|\gamma|\leq T}\frac{x^{i\gamma}}{\frac{1}{2}+i\gamma}-\sum_{|\gamma|\leq T_0}\frac{x^{i\gamma}}{\frac{1}{2}+i\gamma}.
\end{equation}
Then we have the following result (See Section \ref{subsec-lower-trunc} for the proof).
\begin{lem}\label{lem-distri-extra}
	For the lower order sum  and error  from the truncation, we have
	\begin{equation*}
	\int_2^Y \left| S_1(e^y; \chi) \right|^2 dy=o\( Y(\log Y)^{2k-2}\),
	\end{equation*}
	and	for fixed large $T_0$,
	\begin{equation*}
	\int_2^Y |S_2(e^y, T_0; \chi)|^2 dy\ll Y \frac{\log^2 T_0}{T_0} +\log Y \frac{\log^3 T_0}{T_0}+\log^5 T_0.
	\end{equation*}
\end{lem}

Combining Lemmas \ref{lem-bias}, \ref{lem-distri-sum}, and \ref{lem-distri-extra} with (\ref{psi-chi-Omega-form}), (\ref{psi-chi-omega-form}), (\ref{sum-rho-form}), and (\ref{sum-rho-form-omega}), we get, for fixed large $T_0$,
\begin{align}\label{Psi-Omega}
\psi_{\Omega_k}(x, \chi)&=\frac{(-1)^k}{(k-1)!} \frac{\sqrt{x}}{\log x} (\log\log x)^{k-1}\(\sum_{|\gamma|\leq T_0}\frac{x^{i\gamma}}{\frac{1}{2}+i\gamma}+\Sigma(x, T_0; \chi)\) \nonumber\\
&\quad + a(\chi)\frac{(-1)^k}{(k-1)!} \frac{\sqrt{x}}{\log x} (\log\log x)^{k-1},
\end{align}
where
$$\limsup_{Y\rightarrow\infty} \frac{1}{Y} \int_1^Y \left|\Sigma(e^y, T_0; \chi) \right|^2 dy\ll \frac{\log^2 T_0}{T_0}.$$

Also,
\begin{align}\label{Psi-omega}
\psi_{\omega_k}(x, \chi)&=\frac{(-1)^k}{(k-1)!} \frac{\sqrt{x}}{\log x} (\log\log x)^{k-1}\(\sum_{|\gamma|\leq T_0}\frac{x^{i\gamma}}{\frac{1}{2}+i\gamma}+\widetilde{\Sigma}(x, T_0; \chi)\) \nonumber\\
&\quad + a(\chi)\frac{1}{(k-1)!} \frac{\sqrt{x}}{\log x} (\log\log x)^{k-1},
\end{align}
where
$$\limsup_{Y\rightarrow\infty} \frac{1}{Y} \int_1^Y \left|\widetilde{\Sigma}(e^y, T_0; \chi) \right|^2 dy\ll\frac{\log^2 T_0}{T_0}.$$

Note that
$\sum_{\chi\neq \chi_0} (\overline{\chi}(a)-\overline{\chi}(b))a(\chi)=N(q, a)-N(q, b).$
Hence, combining (\ref{Psi-Omega}) and (\ref{Psi-omega}) with (\ref{Delta-Omega}) and (\ref{Delta-omega}), we get the conclusions of Theorem \ref{thm-Delta-Omega} and Theorem \ref{thm-Delta-omega}. \qed

\section{The bias terms}\label{subsec-bias}
In this section, we examine the bias terms and give the proof of Lemma \ref{lem-bias}.

\subsection{Estimates  on the horizontal line}
In order to examine the corresponding integration on the horizontal line in the Hankel contour, we prove the following estimate which we will use many times later to analyze the error terms in our theorems.

\begin{lem}\label{lem-Esigma}
	For any integers $k\geq 1$ and $m\geq 0$, we have
	\begin{equation*}
	\int_0^{\delta} |(\log \sigma -i\pi)^k-(\log\sigma +i\pi)^k|\sigma^m x^{-\sigma}d\sigma\ll_{m,k} \frac{(\log\log x)^{k-1}}{(\log x)^{m+1}}.
	\end{equation*}
\end{lem}
\noindent \textit{\textbf{Proof. }} Let $I$ represent the integral in the lemma. Then, we have
\begin{equation}\label{pf-lem-Esigma-I}
I\leq 2\sum_{j=1}^k {k\choose j} \pi^j \int_0^{\delta} |\log\sigma|^{k-j}\sigma^m x^{-\sigma}d\sigma\ll_k \sum_{j=1}^k \int_0^{\delta} |\log\sigma|^{k-j}\sigma^m x^{-\sigma}d\sigma
\end{equation}
Using a change of variable, $\sigma\log x=t$, we have
\begin{align}\label{pf-lem-Esigma-I-j-for}
&\int_0^{\delta} |\log\sigma|^{k-j}\sigma^m x^{-\sigma}d\sigma\leq \frac{1}{(\log x)^{m+1}} \int_0^{\delta\log x}|\log t-\log\log x|^{k-j}t^m e^{-t}dt\nonumber\\
&\quad \leq \frac{1}{(\log x)^{m+1}} \sum_{l=0}^{k-j} {k-j \choose l} (\log\log x)^{k-j-l} \int_0^{\delta\log x} |\log t|^l t^{m} e^{-t}dt\nonumber\\
&\quad \ll_k  \frac{1}{(\log x)^{m+1}} \sum_{l=0}^{k-j} (\log\log x)^{k-j-l} \int_0^{\delta\log x} |\log t|^l t^{m} e^{-t}dt.
\end{align}
Next, we estimate
\begin{equation}\label{pf-lem-Esigma-I-j-l}
\int_0^{\delta\log x} |\log t|^l t^m e^{-t}dt\leq \(\int_0^1+\int_1^{\infty}  \) |\log t|^l t^m e^{-t} dt=: I_{l_1}+I_{l_2}.
\end{equation}
For the first integral in (\ref{pf-lem-Esigma-I-j-l}),
\begin{equation*}
I_{l_1}=\int_0^1 |\log t|^l t^m e^{-t} dt\leq \int_0^1 |\log t|^l dt\stackrel{t\rightarrow \frac{1}{e^t}}{=}
\int_0^{\infty} \frac{t^l}{e^t}dt =\Gamma(l+1).
\end{equation*}
For the second integral in (\ref{pf-lem-Esigma-I-j-l}),
\begin{equation}\label{pf-lem-Esigma-Il2}
I_{l_2}=\int_1^{\infty} \frac{t^m (\log t)^l}{e^t} dt \stackrel{t\rightarrow e^t}{=} \int_0^{\infty} \frac{t^l}{e^{e^t-(m+1)t}} dt \ll_m \Gamma(l+1).
\end{equation}
Then, by (\ref{pf-lem-Esigma-I-j-for})-(\ref{pf-lem-Esigma-Il2}), we have
\begin{equation}\label{pf-lem-Esigma-I-j}
\int_0^{\delta} |\log\sigma|^{k-j}\sigma^m x^{-\sigma}d\sigma\ll_k  \frac{1}{(\log x)^{m+1}} \sum_{l=0}^{k-j} (\log\log x)^{k-j-l}O_{m,l}(1)\ll_{m, k}\frac{(\log\log x)^{k-j}}{(\log x)^{m+1}}.
\end{equation}
Thus, by (\ref{pf-lem-Esigma-I}),
\begin{eqnarray*}
	I\ll_{m,k} \frac{(\log\log x)^{k-1}}{(\log x)^{m+1}}.
\end{eqnarray*}
Hence, we get the conclusion of this lemma. \qed

\subsection{The bias terms}

We have the following estimate for the integral over the truncated Hankel contour $\mathcal{H}(\frac{1}{2}, \delta)$. 
\begin{lem} \label{lem-bias-err}
	Assume the function $f(s)=O(1)$ on $\mathcal{H}(\frac{1}{2},\delta)$. Then, for any integer $m\geq 0$,
	\begin{equation*}
	\left|\int_{\mathcal{H}(\frac{1}{2},\delta)} \(\log \(s-\frac{1}{2} \) \)^{m} f(s) \frac{x^s}{s}ds\right| \ll_m \frac{\sqrt{x}(\log\log x)^{m-1}}{\log x}.
	\end{equation*}
\end{lem}
\noindent \textbf{\textit{Proof.}} Since the left-hand side is $0$ when $m=0$, we assume $m\geq 1$ in the following proof. By Lemma \ref{lem-Esigma}, we have
\begin{align}\label{pf-bias-err-1-j}
&\left|\int_{\mathcal{H}(\frac{1}{2},\delta)} \(\log \(s-\frac{1}{2} \) \)^{m} f(s) \frac{x^s}{s}ds\right|\nonumber\\
& \leq\left|\int_{r_0}^{\delta} \((\log\sigma-i\pi)^{m}-(\log\sigma+i\pi)^{m} \) f\(\frac{1}{2}-\sigma\) \frac{x^{\frac{1}{2}-\sigma}}{\frac{1}{2}-\sigma}d\sigma\right|\nonumber\\
&\quad +O \(\int_{-\pi}^{\pi}  \frac{\(\log\frac{1}{r_0}+\pi\)^{m}x^{\frac{1}{2}+r_0} }{\frac{1}{2}-r_0} r_0 d\alpha \)\nonumber\\
& \ll \sqrt{x} \(\int_0^{\delta} \left|(\log\sigma-i\pi)^{m}-(\log\sigma+i\pi)^{m} \right| x^{-\sigma}d\sigma +\frac{(\log x+\pi)^{m}}{x}\)\nonumber\\
&\ll_m  \frac{\sqrt{x}(\log\log x)^{m-1}}{\log x}.
\end{align}
This completes the proof of this lemma. \qed

\bigskip

In the following, we prove the asymptotic formulas for the bias terms. 

\noindent \textit{\textbf{Proof of Lemma \ref{lem-bias}.}}
Since $H_B(s)=O(1)$, by (\ref{bias-form}), (\ref{bias-form-omega}), and Lemma \ref{lem-bias-err},
\begin{align}\label{bias-err-1}
|E_B(x)|&\ll\sum_{j=1}^k\left|  \int_{\mathcal{H}(\frac{1}{2},\delta)} \(\log \(s-\frac{1}{2} \) \)^{k-j} \(H_B(s) \)^j \frac{x^s}{s}ds\right|\nonumber\\
&\ll \frac{\sqrt{x}}{\log x} \sum_{j=1}^k  (\log\log x)^{k-j-1}\ll_k \frac{\sqrt{x}}{\log x} (\log\log x)^{k-2}.
\end{align}
Similarly,
\begin{equation}\label{bias-err-omega-1}
|\widetilde{E}_B(x)|\ll\sum_{j=1}^k\left|  \int_{\mathcal{H}(\frac{1}{2},\delta)} \(\log \(s-\frac{1}{2} \) \)^{k-j} \(\widetilde{H}_B(s) \)^j \frac{x^s}{s}ds\right|\ll_k \frac{\sqrt{x}}{\log x} (\log\log x)^{k-2}.
\end{equation}

In the following, we estimate $E_R(x)$ in (\ref{bias-form}) and $\widetilde{E}_R(x)$ in (\ref{bias-form-omega}). If $\chi$ is not real, $E_R(x)=\widetilde{E}_R(x)=0$. If $\chi$ is real, by (\ref{F-to-L}), on $\mathcal{H}(\frac{1}{2}, \delta)$, we write
\begin{equation}\label{bias-F2s-form}
F(2s, \chi^2)=-\log\(s-\frac{1}{2}\)+H_2(s).
\end{equation}
On $\mathcal{H}(\frac{1}{2}, \delta)$, $|H_2(s)|=O(1)$. By (\ref{F-Omega-k}), we have
\begin{equation}\label{pf-bias-ER-sum}
|E_{R}(x)|\ll \sum_{m=0}^{k-2} \sum_{\boldsymbol{n}\in S^{(k)}_m}\left|\int_{\mathcal{H}(\frac{1}{2}, \delta)} F^m(s, \chi) F(\boldsymbol{n}s, \chi)\frac{x^s}{s}ds  \right|.
\end{equation}
For each $0\leq m\leq k-2$, we write
\begin{equation*}
F^m(s, \chi) F(\boldsymbol{n}s, \chi)=F^m(s, \chi) F^{m'}(2s,\chi^2)G_{\boldsymbol{n}}(s),
\end{equation*}
where $m+2m'\leq k$, and $G_{\boldsymbol{n}}(s)=O(1)$ on $\mathcal{H}(\frac{1}{2}, \delta)$.
Thus, by (\ref{pf-thm-bias-F-formu}), (\ref{bias-F2s-form}), and Lemma \ref{lem-bias-err},
\begin{eqnarray}\label{pf-bias-ER-con}
&&\left|\int_{\mathcal{H}(\frac{1}{2}, \delta)} F^m(s, \chi) F(\boldsymbol{n}s, \chi)\frac{x^s}{s}ds  \right| \nonumber\\
&& \ll \left|  \int_{\mathcal{H}(\frac{1}{2},\delta)} \(\log\(s-\frac{1}{2} \) +H_B(s)\)^{m} \(\log\(s-\frac{1}{2} \)-H_2(s) \)^{m'}G_{\boldsymbol{n}}(s)\frac{x^s}{s}ds\right|\nonumber\\
&&  \ll \sum_{j_1=0}^m \sum_{j_2=0}^{m'} \left|\int_{\mathcal{H}(\frac{1}{2}, \delta)}\(\log\(s-\frac{1}{2}\) \)^{m+m'-j_1-j_2} (H_B(s))^{j_1} (H_2(s))^{j_2} G_{\boldsymbol{n}}(s) \frac{x^s}{s}ds \right|\nonumber\\
&& \ll \sum_{j_1=0}^m \sum_{j_2=0}^{m'} \frac{\sqrt{x}}{\log x} (\log\log x)^{m+m'-j_1-j_2-1} \ll_k \frac{\sqrt{x}}{\log x} (\log\log x)^{k-2}.
\end{eqnarray}
In the last step, we used the conditions $0\leq m\leq k-2$ and $m+2m'\leq k$.

Combining  (\ref{pf-bias-ER-sum}) and (\ref{pf-bias-ER-con}),  we deduce that
\begin{equation}\label{bias-ER}
|E_R(x)|\ll_k \frac{ \sqrt{x}}{\log x} (\log\log x)^{k-2}.
\end{equation}

Similarly, if $\chi$ is real, by (\ref{F-omega-L-2s}), we write
\begin{equation}\label{bias-F2s-form-omega}
\widetilde{F}(s, \chi; 2)=-\log\(s-\frac{1}{2}\)+\widetilde{H}_2(s),
\end{equation}
where $\widetilde{H}_2(s)=O(1)$ on $\mathcal{H}(\frac{1}{2}, \delta)$. Using a similar argument as above, by (\ref{pf-thm-bias-F-omega-formu}), (\ref{bias-F2s-form-omega}), and Lemma \ref{lem-bias-err}, we have
\begin{equation}\label{bias-ER-omega}
|\widetilde{E}_R(x)|\ll \sum_{m=0}^{k-2} \sum_{\boldsymbol{n}\in S^{(k)}_m}\left|\int_{\mathcal{H}(\frac{1}{2}, \delta)} \widetilde{F}^m(s, \chi) \widetilde{F}(\boldsymbol{n}s, \chi)\frac{x^s}{s}ds  \right|
\ll_k \frac{ \sqrt{x}}{\log x} (\log\log x)^{k-2}.
\end{equation}

By (\ref{bias-form}), (\ref{main-bias}), (\ref{bias-err-1}), and (\ref{bias-ER}), we get
\begin{equation*}
I_B(x)=\frac{(-1)^k\sqrt{x}}{2^{k-1}\log x}\left\{  k (\log\log x)^ {k-1}+\sum_{j=2}^k {k \choose j} (\log\log x)^{k-j} \frac{1}{\Gamma_j(0)}\right\}+O_k\(\frac{\sqrt{x}(\log\log x)^{k-2}}{\log x}\).
\end{equation*}
Then, by (\ref{Gamma-j-estimate}),
$$\left|\sum_{j=2}^k {k \choose j} (\log\log x)^{k-j} \frac{1}{\Gamma_j(0)}\right|\ll_k (\log\log x)^{k-2}.$$
Hence,
\begin{equation}\label{bias-IB}
I_B(x)=\frac{(-1)^k k}{2^{k-1}}\frac{\sqrt{x}}{\log x} (\log\log x)^{k-1}+O_k\( \frac{\sqrt{x}(\log\log x)^{k-2}}{\log x}\).
\end{equation}

Similarly, by (\ref{bias-form-omega}), (\ref{main-bias-omega}), (\ref{bias-err-omega-1}), and (\ref{bias-ER-omega}), we have
\begin{equation}\label{bias-IB-omega}
\widetilde{I}_B(x)=\frac{ k}{2^{k-1}}\frac{\sqrt{x}}{\log x} (\log\log x)^{k-1}+O_k\( \frac{\sqrt{x}(\log\log x)^{k-2}}{\log x}\).
\end{equation}

This completes the proof of Lemma \ref{lem-bias}. \qed

\section{Average order of the error terms}\label{sec-error-terms}
In Section \ref{subsec-error-Hankel} and Section \ref{subsec-esti-Hankel}, we examine the error terms from the Hankel contours around zeros and give the proof of Lemma \ref{lem-distri-sum}. In Section \ref{subsec-lower-trunc}, we examine the lower order sum and the error from the truncation, and give the proof of Lemma \ref{lem-distri-extra}. 

\subsection{Error terms from the Hankel contours around zeros}\label{subsec-error-Hankel}
In the section, we give the proof of Lemma \ref{lem-distri-sum}. The following lemma gives an average estimate for the integral over Hankel contours around zeros, which is the key lemma for our proof. 
 \begin{lem}\label{pf-distri-Err-gene-lem}
 	Let $\rho$ be a zero of $L(s, \chi)$. Assume the function $g(s)\ll \( \log |\gamma|\)^c$ on $\mathcal{H}(\rho, \delta)$ for some constant $c\geq 0$, and
 	\begin{equation}\label{lem-distri-sum-func-H}
 	H_{\rho}(s)=\sum_{0<|\gamma'-\gamma|\leq 1} \log (s-\rho') +O\(\log |\gamma| \)\quad{\rm on~} \mathcal{H}(\rho,\delta).
 	\end{equation}
 	For any integers $m, n\geq 0$, denote
 	\begin{equation*}
 	E(x; \rho):=\int_{\mathcal{H}(\rho,\delta)} \(\log(s-\rho)\)^{m} \(H_{\rho}(s)\)^n g(s) \frac{x^{s-\rho}}{s}ds.
 	\end{equation*}
 	Then, for $T=T(x)$, we have
 	\begin{equation*}
 	\int_{2}^Y \left|y \sum_{|\gamma|\leq T(e^y)} e^{i\gamma y} E(e^y; \rho) \right|^2 dy=o\( Y(\log Y)^{2m+2n-2} \).
 	\end{equation*}
 \end{lem}
We will give the proof of Lemma \ref{pf-distri-Err-gene-lem} in next subsection.  We use it in this section to prove Lemma \ref{lem-distri-sum} first.

\vspace{1em}
\noindent\textit{\textbf{Proof of Lemma \ref{lem-distri-sum}.}}
By (\ref{de-Sigma-1}), we have
\begin{equation}\label{pf-distri-EM-formu}
\left| \Sigma_1(x; \chi) \right|^2 =\left|\log x\sum_{|\gamma|\leq T} x^{i\gamma} E'_{M_{\rho}}(x)\right|^2
\ll \sum_{j=1}^k \left| \log x \sum_{|\gamma|\leq T} x^{i\gamma} E_{\rho, j}(x)  \right|^2,
\end{equation}
where
$$E_{\rho, j}(x)=\int_{\mathcal{H}(\rho,\delta)} (\log(s-\rho))^{k-j} (H_{\rho}(s))^j \frac{x^{s-\rho}}{s}ds.$$
By Lemma \ref{pf-distri-Err-gene-lem}, take $m=k-j$, $n=j$, and $g(s)\equiv 1$, (i.e. $c=0$),
\begin{equation*}
\int_2^Y \left|y \sum_{|\gamma|\leq T(e^y)} e^{i\gamma y} E_{\rho, j}(e^y)  \right|^2 dy=o\(Y(\log Y)^{2k-2}\).
\end{equation*}
Thus,
\begin{equation}\label{pf-distri-EM-Omega-Y}
\int_2^Y \left| \Sigma_1(e^y; \chi) \right|^2 dy=o\(Y(\log Y)^{2k-2}\).
\end{equation}

By definition (\ref{rho-form}) and (\ref{de-Sigma-2}), we have
\begin{align}\label{pf-distri-ER-formu}
\left| \Sigma_2(x; \chi) \right|^2&\ll\sum_{m=0}^{k-2} \sum_{\boldsymbol{n}\in S^{(k)}_m}\left| \log x \sum_{|\gamma|\leq T} x^{i\gamma}\int_{\mathcal{H}(\rho, \delta)} F^m(s, \chi)F(\boldsymbol{n}s, \chi) \frac{x^{s-\rho}}{s}ds \right|^2\nonumber\\
&\ll\sum_{m=0}^{k-2} \sum_{\boldsymbol{n}\in S^{(k)}_m} \sum_{j=0}^m \left|\log x \sum_{|\gamma|\leq T} x^{i\gamma} E_{m, j}(x, \chi; \boldsymbol{n}) \right|^2,
\end{align}
where
$$E_{m,j}(x, \chi; \boldsymbol{n})=\int_{\mathcal{H}(\rho, \delta)} \(\log(s-\rho) \)^{m-j} \(H_{\rho}(s) \)^j F(\boldsymbol{n}s, \chi) \frac{x^{s-\rho}}{s}ds.$$
Since on $\mathcal{H}(\rho, \delta)$, we know $F(\boldsymbol{n}s, \chi)=O\((\log |\gamma|)^{\frac{k-m}{2}} \)$ , by Lemma \ref{pf-distri-Err-gene-lem}, we get
\begin{equation*}
\int_2^Y \left|y \sum_{|\gamma|\leq T(e^y)} e^{i\gamma y} E_{m, j}(e^y, \chi; \boldsymbol{n}) \right|^2 dy=o\( Y(\log Y)^{2m-2}\).
\end{equation*}
Hence, by (\ref{pf-distri-ER-formu}), we deduce that
\begin{equation}\label{pf-distri-ER-Omega-Y}
\int_2^Y \left| \Sigma_2(e^y; \chi) \right|^2 dy =o\(Y(\log Y)^{2k-2} \).
\end{equation}
Combining (\ref{pf-distri-EM-Omega-Y}) and (\ref{pf-distri-ER-Omega-Y}), we get the first formula in Lemma \ref{lem-distri-sum}.

For $\widetilde{\Sigma}_1(x; \chi)$ and $\widetilde{\Sigma}_2(x, \chi)$, by (\ref{de-Sigma-1-omega}), using a similar argument with Lemma \ref{pf-distri-Err-gene-lem},
\begin{equation}\label{pf-distri-EM-formu-omega}
\int_2^Y \left| \widetilde{\Sigma}_1(e^y; \chi) \right|^2 dy
\ll \sum_{j=1}^k \int_2^Y \left| y \sum_{|\gamma|\leq T(e^y)} e^{i\gamma y} \widetilde{E}_{\rho, j}(e^y)  \right|^2 dy=o\(Y(\log Y)^{2k-2} \),
\end{equation}
where
$$\widetilde{E}_{\rho, j}(x)=\int_{\mathcal{H}(\rho,\delta)} (\log(s-\rho))^{k-j} (\widetilde{H}_{\rho}(s))^j \frac{x^{s-\rho}}{s}ds.$$

Similarly, by (\ref{de-Sigma-2-omega}) and Lemma \ref{pf-distri-Err-gene-lem},
\begin{equation}\label{pf-distri-ER-formu-omega}
\int_2^Y \left| \widetilde{\Sigma}_2(e^y; \chi) \right|^2 dy
\ll\sum_{m=0}^{k-2} \sum_{\boldsymbol{n}\in S^{(k)}_m} \sum_{j=0}^m \int_2^Y \left|y \sum_{|\gamma|\leq T(e^y)} e^{i\gamma y} \widetilde{E}_{m, j}(e^y, \chi; \boldsymbol{n}) \right|^2 dy= o\(Y(\log Y)^{2k-2} \),
\end{equation}
where
$$\widetilde{E}_{m,j}(x, \chi; \boldsymbol{n})=\int_{\mathcal{H}(\rho, \delta)} \(\log(s-\rho) \)^{m-j} \(\widetilde{H}_{\rho}(s) \)^j \widetilde{F}(\boldsymbol{n}s, \chi) \frac{x^{s-\rho}}{s}ds.$$

Combining (\ref{pf-distri-EM-formu-omega}) and (\ref{pf-distri-ER-formu-omega}), we get the second formula in Lemma \ref{lem-distri-sum}. \qed

\subsection{Estimates for integral over Hankel contours around zeros}\label{subsec-esti-Hankel}
We  need the following results to finish the proof of Lemma \ref{pf-distri-Err-gene-lem}. 

\begin{lem}[\cite{Ford-S}, Lemma 2.4] \label{lem-Ford-S}
	Assume $L(\frac{1}{2}, \chi)\neq 0$. For $A\geq 0$ and real $l\geq 0$,
	\begin{equation*}
	\sum_{\substack{|\gamma_1|, |\gamma_2|\geq A\\ |\gamma_1-\gamma_2|\geq 1}} \frac{\log^l (|\gamma_1|+3) \log^l (|\gamma_2| +3)}{|\gamma_1||\gamma_2||\gamma_1-\gamma_2|}\ll_l \frac{(\log(A+3))^{2l+3}}{A+1}.
	\end{equation*}
\end{lem}

\begin{lem}\label{lem-Delta-N}
	For any integers $N, j \geq 1$, and $0<|\delta_n|\leq 1$, we have
	\begin{equation*}
	\int_0^{\delta}\left| \sum_{n=1}^N \log(\sigma+i\delta_n)\right|^j x^{-\sigma}d\sigma \ll_j \frac{1}{\log x}\left\{\min\(N\log\log x, \log\frac{1}{\Delta_N}\)+N\pi \right\}^j,
	\end{equation*}
	where $\Delta_N=\prod_{n=1}^N |\delta_n|$.
\end{lem}
\noindent\textit{\textbf{Proof.}} Let $I$ denote the integral in the lemma. We consider two cases: $\Delta_N\geq \(\frac{1}{\log x}\)^N$, and $\Delta_N< \(\frac{1}{\log x}\)^N$.

1) If $\Delta_N\geq \(\frac{1}{\log x}\)^N$, we have
\begin{equation}\label{pf-lem-E-Delta-N-case1-c}
I\ll \(\log\frac{1}{\Delta_N}+N\pi\)^j \int_0^{\delta} x^{-\sigma}d\sigma\ll \frac{1}{\log x} \(\log\frac{1}{\Delta_N}+N\pi\)^j
\ll \frac{1}{\log x}(N\log\log x+N\pi)^j.
\end{equation}

2) If $\Delta_N< \(\frac{1}{\log x}\)^N$, we write
\begin{equation}\label{pf-lem-E-Delta-N-case2-I}
I=\(\int_0^{(\Delta_N)^{\frac{1}{N}}}+\int_{(\Delta_N)^{\frac{1}{N}}}^{\frac{1}{\log x}}+\int_{\frac{1}{\log x}}^{\delta} \)\left| \sum_{n=1}^N \log(\sigma+i\delta_n)\right|^j x^{-\sigma}d\sigma
=: I_1+I_2+I_3.
\end{equation}
First, we estimate $I_1$,
\begin{equation}\label{pf-lem-E-Delta-N-I1}
I_1\ll\(\log\frac{1}{\Delta_N}+N\pi\)^j \int_0^{(\Delta_N)^{\frac{1}{N}}} x^{-\sigma}d\sigma
\ll (\Delta_N)^{\frac{1}{N}} \(\log\frac{1}{\Delta_N}+N\pi\)^j.
\end{equation}
For $0<t<1$, consider the function
$f(t)=t^{\frac{1}{N}}\(\log\frac{1}{t}+N\pi\)^j.$
Since the critical point of $f(t)$ is $t=e^{N(\pi-1)}>1$, by (\ref{pf-lem-E-Delta-N-I1}), we have
\begin{equation}\label{pf-lem-Delta-N-I1-c}
I_1\ll f\(\frac{1}{(\log x)^N}\)=\frac{1}{\log x}(N\log\log x+N\pi)^j\ll \frac{1}{\log x} \(\log\frac{1}{\Delta_N}+N\pi\)^j.
\end{equation}

Next, we estimate $I_3$. Using the change of variable $\sigma\log x=t$, we get
\begin{eqnarray}\label{pf-lem-E-Delta-N-I3}
I_3&\ll& \int_{\frac{1}{\log x}}^{\delta} \(N\log\frac{1}{\sigma}+N\pi\)^j x^{-\sigma} d\sigma \nonumber\\
&=& \frac{1}{\log x}\int_1^{\delta\log x} (N\log\log x-N\log t+N\pi)^j e^{-t} dt\nonumber\\
&=& \frac{N^j}{\log x} \sum_{l=0}^j {j\choose l} (\log\log x+\pi)^{j-l} \int_1^{\delta\log x} (-\log t)^l e^{-t} dt\nonumber\\
&\ll_j& \frac{N^j}{\log x} \sum_{l=0}^j  (\log\log x+\pi)^{j-l} \int_1^{\infty}\frac{t^l}{e^t}dt\nonumber\\
&\ll_j& \frac{ (N\log\log x+N\pi)^j}{\log x}\ll \frac{1}{\log x} \(\log\frac{1}{\Delta_N}+N\pi\)^j.
\end{eqnarray}

For $I_2$, similar to $I_3$, using the change of variable $\sigma\log x=t$, we get
\begin{eqnarray}\label{pf-lem-E-Delta-N-I2}
I_2&\ll& \int_{(\Delta_N)^{\frac{1}{N}}}^{\frac{1}{\log x}} \(N\log\frac{1}{\sigma}+N\pi\)^j x^{-\sigma} d\sigma \nonumber\\
&=& \frac{1}{\log x}\int_{(\Delta_N)^{\frac{1}{N}}\log x}^1 (N\log\log x-N\log t+N\pi)^j e^{-t} dt\nonumber\\
&=& \frac{N^j}{\log x} \sum_{l=0}^j {j\choose l} (\log\log x+\pi)^{j-l} \int_{(\Delta_N)^{\frac{1}{N}}\log x}^1 (-\log t)^l e^{-t} dt \quad (t\rightarrow \frac{1}{e^t}) \nonumber\\
&\ll_j& \frac{N^j}{\log x} \sum_{l=0}^j  (\log\log x+\pi)^{j-l} \int_0^{\infty} \frac{t^l}{e^t}dt\nonumber\\
&\ll_j& \frac{ (N\log\log x+N\pi)^j}{\log x} \ll \frac{1}{\log x} \(\log\frac{1}{\Delta_N}+N\pi\)^j.
\end{eqnarray}
Combining (\ref{pf-lem-Delta-N-I1-c}), (\ref{pf-lem-E-Delta-N-I3}), (\ref{pf-lem-E-Delta-N-I2}), with (\ref{pf-lem-E-Delta-N-case2-I}), we get
\begin{equation}\label{pf-lem-E-Delta-N-case2-c}
I\ll_j \frac{ (N\log\log x+N\pi)^j}{\log x} \ll_j \frac{1}{\log x} \(\log\frac{1}{\Delta_N}+N\pi\)^j.
\end{equation}

By (\ref{pf-lem-E-Delta-N-case1-c}) and (\ref{pf-lem-E-Delta-N-case2-c}), we get the conclusion of this lemma. \qed

\bigskip
In the following, we use the above lemmas to prove Lemma \ref{pf-distri-Err-gene-lem}. 

 \noindent \textit{\textbf{Proof of Lemma \ref{pf-distri-Err-gene-lem}.}}
 If $m=0$, $E(x; \rho)=0$ and hence the integral is $0$. In the following, we assume $m\geq 1$.
 Let $\Gamma_{\rho}$ represent the circle in the Hankel contour $\mathcal{H}(\rho, \delta)$. Then,
 \begin{align}\label{pf-distri-Err-lem-E-formu}
 E(x; \rho)&=\int_{\mathcal{H}(\rho,\delta)} \(\log(s-\rho)\)^{m} \(H_{\rho}(s)\)^n g(s) \frac{x^{s-\rho}}{s}ds  \nonumber\\
 &= \int_{r_{\rho}}^{\delta} \((\log\sigma-i\pi)^{m}-(\log\sigma+i\pi)^{m} \) \(H_{\rho}\(\frac{1}{2}-\sigma+i\gamma\)\)^n g\(\frac{1}{2}-\sigma+i\gamma\)\nonumber\\
 & \quad\quad \times \frac{x^{-\sigma}}{\frac{1}{2}-\sigma+i\gamma}d\sigma+\int_{\Gamma_{\rho}} \(\log(s-\rho)\)^{m} \(H_{\rho}(s)\)^n g(s) \frac{x^{s-\rho}}{s}ds.\nonumber\\
 &=: E_h(x; \rho)+E_r(x; \rho).
 \end{align}
 For the second integral in (\ref{pf-distri-Err-lem-E-formu}), since  $r_{\rho}\leq \frac{1}{x}$, by Lemma \ref{lem-Ford-S-2.2},
 \begin{align}\label{pf-distri-Err-lem-E-r-est}
 \left|E_r(x; \rho)  \right|
 &\ll \frac{(\log |\gamma|)^c r_{\rho} x^{r_{\rho}} }{|\gamma|} \(\log\frac{1}{r_{\rho}} +\pi\)^{m} \(\sum_{0<|\gamma-\gamma'|\leq 1} \log\(\frac{1}{|\gamma'-\gamma|-r_{\rho}}\) +O(\log|\gamma|) \)^n\nonumber\\
 &\ll \frac{(\log |\gamma|)^c r_{\rho} x^{r_{\rho}} }{|\gamma|} \(\log\frac{1}{r_{\rho}} +\pi\)^{m} (\log|\gamma|)^n\( \log\(\frac{1}{r_{\rho}} \)+O(1) \)^n\nonumber\\
 &\ll\frac{(\log|\gamma| )^{n+c}}{|\gamma|} \frac{  \(\log(1/r_{\rho})+\pi\)^{m+n}}{1/r_{\rho}}\ll \frac{(\log|\gamma| )^{n+c}}{|\gamma|} \frac{1}{x^{1-\epsilon}}.
 \end{align}

Denote
\begin{equation}\label{pf-distri-Err-lem-Sigma-h-r}
\Sigma(x; {\rm g}):= \left| \sum_{|\gamma|\leq T} x^{i\gamma} E(x; \rho) \right|^2
\ll\left| \sum_{|\gamma|\leq T} x^{i\gamma} E_h(x; \rho) \right|^2+\left| \sum_{|\gamma|\leq T} x^{i\gamma} E_r(x; \rho) \right|^2.
\end{equation}
By (\ref{pf-distri-Err-lem-E-r-est}),  and $ T(x)\ll x^2$, we get
\begin{equation}\label{pf-distri-Err-lem-Sigma-r-est}
\left| \sum_{|\gamma|\leq T} x^{i\gamma} E_r(x; \rho) \right|^2\ll\frac{1}{x^{2-\epsilon}}\(\sum_{|\gamma|\leq T(x)} \frac{(\log|\gamma| )^{n+c}}{|\gamma|}\)^2  \ll \frac{1}{x^{2-\epsilon}}.
\end{equation}
For the first sum in (\ref{pf-distri-Err-lem-Sigma-h-r}),
 \begin{align*}
 \left| \sum_{|\gamma|\leq T} x^{i\gamma} E_h(x; \rho) \right|^2
&=\(\sum_{\substack{|\gamma_1-\gamma_2|\leq 1\\|\gamma_1|, |\gamma_2|\leq T}}+\sum_{\substack{|\gamma_1-\gamma_2|> 1\\|\gamma_1|, |\gamma_2|\leq T}}\)  x^{i(\gamma_1-\gamma_2)} E_h(x; \rho_1) E_h(x; \overline{\rho}_2)\nonumber\\
 &=:\Sigma_1(x; {\rm g})+\Sigma_2(x; {\rm g}).
 \end{align*}


By (\ref{pf-distri-Err-lem-E-formu}),
 \begin{equation}\label{pf-distri-Err-lem-E-h}
 \left| E_h(x; \rho)\right| \ll \frac{(\log |\gamma|)^c}{|\gamma|} \sum_{j=1}^{m}  \int_0^{\delta} |\log \sigma|^{m-j} \left|H_{\rho}\(\frac{1}{2}-\sigma+i\gamma\) \right|^n x^{-\sigma}d\sigma.
 \end{equation}
 Let
 \begin{align*}
 S_j(x)&:=\int_0^{\delta} |\log \sigma|^{m-j} \left|H_{\rho}\(\frac{1}{2}-\sigma+i\gamma\) \right|^n x^{-\sigma}d\sigma\nonumber\\
 &\leq \(\int_0^{\delta} |\log \sigma|^{2(m-j)} x^{-\sigma}d\sigma\)^{\frac{1}{2}} \(\int_0^{\delta} \left|H_{\rho}\(\frac{1}{2}-\sigma+i\gamma\) \right|^{2n} x^{-\sigma}d\sigma \)^{\frac{1}{2}}.
 \end{align*}
 By (\ref{pf-lem-Esigma-I-j}) in the proof of Lemma \ref{lem-Esigma},
 \begin{equation}\label{pf-distri-Er-SL-1}
 \int_0^{\delta} |\log \sigma|^{2(m-j)} x^{-\sigma}d\sigma\ll \frac{(\log\log x)^{2(m-j)}}{\log x}.
 \end{equation}
 By condition (\ref{lem-distri-sum-func-H}) and the Cauchy-Schwarz inequality,
 \begin{equation*}
 \left|H_{\rho}\(\frac{1}{2}-\sigma+i\gamma\) \right|^{2n}\ll \left|\sum_{0<|\gamma'-\gamma|\leq 1} \log(\sigma+i(\gamma'-\gamma))\right|^{2n}+(\log|\gamma|)^{2n}.
 \end{equation*}
 Then, by Lemma \ref{lem-Delta-N},
 \begin{equation}\label{pf-distri-Er-SL-H}
 \int_0^{\delta} \left|H_{\rho}\(\frac{1}{2}-\sigma+i\gamma\) \right|^{2n} x^{-\sigma}d\sigma \ll \frac{(M_{\gamma}(x))^{2n}+(\log|\gamma|)^{2n}}{\log x},
 \end{equation}
 where
 $M_{\gamma}(x)=\min\(N(\gamma)\log\log x, \log\frac{1}{\Delta_{N(\gamma)}} \),$
$N(\gamma)$ is the number of zeros $\gamma'$ in the range $ 0<|\gamma'-\gamma|\leq 1 $, and $\Delta_{N(\gamma)}=\prod\limits_{0<|\gamma'-\gamma|\leq 1} |\gamma'-\gamma|$.

 Thus, by (\ref{pf-distri-Er-SL-1}) and (\ref{pf-distri-Er-SL-H}),
 \begin{equation*}
 S_j(x)\ll \frac{(\log\log x)^{m-j}}{\log x}\( (M_{\gamma}(x))^n+(\log|\gamma|)^{n}\).
 \end{equation*}
 Substituting  this into (\ref{pf-distri-Err-lem-E-h}), we get
 \begin{align}\label{pf-distri-Err-lem-E-h-est}
 \left| E_h(x; \rho)  \right|&\ll \frac{(\log |\gamma|)^c}{|\gamma|} \sum_{j=1}^{m}\frac{(\log\log x)^{m-j}}{\log x}\( (M_{\gamma}(x))^n+(\log|\gamma|)^{n}\)\nonumber\\
 &\ll \frac{(\log |\gamma|)^c}{|\gamma|}\frac{(\log\log x)^{m-1}}{\log x}\( (M_{\gamma}(x))^n+(\log|\gamma|)^{n}\).
 \end{align}
Then, by Lemma \ref{lem-Ford-S-2.2}, we have
\begin{align*}
\left|\Sigma_1(x; \rm g)  \right|&\ll \sum_{|\gamma|\leq T} \log(|\gamma|) \(\max_{|\gamma'-\gamma|<1}\left|E_h(x; \rho')\right|\)^2\nonumber\\
&\ll \frac{(\log\log x)^{2(m-1)}}{\log^2 x} \sum_{\gamma}\frac{(\log |\gamma|)^{2c}}{|\gamma|^2}\((M_{\gamma}(x))^{2n}+\(\log|\gamma| \)^{2n} \)\nonumber\\
&= \frac{(\log\log x)^{2m+2n-2}}{\log^2 x} o(1).
\end{align*}
Thus, for each positive integer $l$,
\begin{equation}\label{pf-distri-Err-gene-lem-sum-1}
\int_{2^l}^{2^{l+1}} \Sigma_1(e^y; {\rm g}) dy=o\( \frac{l^{2m+2n-2}}{2^l} \).
\end{equation}


In the following, we examine $\Sigma_2(x; {\rm g})$.  By (\ref{pf-distri-Err-lem-E-formu}),
 \begin{equation}\label{pf-distri-Err-lem-J-sum-formu}
 \Sigma_2(x; {\rm g})=\sum_{\substack{|\gamma_1-\gamma_2|>1\\ |\gamma_1|, |\gamma_2|\leq T}}x^{i(\gamma_1-\gamma_2)} E_h(x; \rho_1)  E_h(x; \overline{\rho}_2).
 \end{equation}
 For $e^{2^l}\leq x\leq e^{2^{l+1}}$, $T=T(x)=T_{l'}$ is a constant, and so we define
 \begin{equation}\label{pf-distri-Err-lem-J-J(x)}
 J(x; {\rm g}):= \sum_{\substack{|\gamma_1-\gamma_2|>1\\|\gamma_1|, |\gamma_2|\leq T_{l'}}} x^{i(\gamma_1-\gamma_2)} \int_{r_{\rho_1}}^{\delta} \int_{r_{\overline{\rho}_2}}^{\delta} R_{\rho_1}(\sigma_1; x) R_{\overline{\rho}_2}(\sigma_2; x) \frac{d\sigma_1 d\sigma_2}{i(\gamma_1-\gamma_2)-(\sigma_1+\sigma_2)},
 \end{equation}
 where
$$ R_{\rho}(\sigma; x)=\((\log\sigma-i\pi)^{m}-(\log\sigma+i\pi)^{m} \) H^n_{\rho}\(\frac{1}{2}-\sigma+i\gamma\) \frac{g\(\frac{1}{2}-\sigma+i\gamma\)x^{-\sigma}}{\frac{1}{2}-\sigma+i\gamma}.$$
Thus,
 \begin{equation}\label{pf-distri-Err-lem-J-sum-to-J}
\int_{e^{2^l}}^{e^{2^{l+1}}} \sum_{\substack{|\gamma_1-\gamma_2|>1\\|\gamma_1|, |\gamma_2|\leq T_{l'}}}  x^{i(\gamma_1-\gamma_2)} E_h(x; \rho_1) E_h(x;\overline{\rho}_2)\frac{dx}{x}=J(e^{2^{l+1}}; {\rm g})-J(e^{2^l}; {\rm g}).
 \end{equation}
By (\ref{pf-distri-Err-lem-J-J(x)}), (\ref{pf-distri-Err-lem-E-h}), and (\ref{pf-distri-Err-lem-E-h-est}), and Lemma \ref{lem-Ford-S}, for $e^{2^l}\leq x\leq e^{2^{l+1}}$
 \begin{align}\label{pf-distri-Err-lem-J-J(x)-est}
 |J(x; {\rm g})|&\ll \sum_{|\gamma_1-\gamma_2|>1} \frac{(\log |\gamma_1|)^c (\log |\gamma_2|)^c}{|\gamma_1||\gamma_2||\gamma_1-\gamma_2|} \(\frac{(\log\log x)^{m-1}}{\log x}\)^2 \nonumber\\
 & \qquad \times \( (M_{\gamma_1}(x))^n+(\log|\gamma_1|)^{n}\)\( (M_{\gamma_2}(x))^n+(\log|\gamma_2|)^{n}\)\nonumber\\
 &\ll \frac{(\log\log x)^{2m+2n-2}}{\log^2 x}\sum_{|\gamma_1-\gamma_2|>1} \frac{(\log |\gamma_1|)^{n+c} (\log |\gamma_2|)^{n+c}}{|\gamma_1||\gamma_2||\gamma_1-\gamma_2|}\ll \frac{(\log\log x)^{2m+2n-2}}{\log^2 x}.
 \end{align}
 Hence, by (\ref{pf-distri-Err-lem-J-sum-formu}), (\ref{pf-distri-Err-lem-J-sum-to-J}), and (\ref{pf-distri-Err-lem-J-J(x)-est}), we get, for any positive integer $l$,
 \begin{equation}\label{pf-distri-Err-gene-lem-sum-2}
 \int_{2^l}^{2^{l+1}} \Sigma_2(e^y; {\rm g}) dy=o\( \frac{l^{2m+2n-2}}{2^l} \).
 \end{equation}

Therefore, by (\ref{pf-distri-Err-lem-Sigma-r-est}), (\ref{pf-distri-Err-gene-lem-sum-1}) and (\ref{pf-distri-Err-gene-lem-sum-2}),
\begin{eqnarray*}
&&\int_{2}^Y \left| y \sum_{|\gamma|\leq T(e^y)} e^{i\gamma y} E(e^y; \rho) \right|^2 dy\ll \sum_{l\leq \frac{\log Y}{\log 2}+1} 2^{2l}
\int_{2^l}^{2^{l+1}} \Sigma(e^y; {\rm g}) dy \\
&&\quad \ll 1+\sum_{l\leq \frac{\log Y}{\log 2}+1} 2^{2l}
\int_{2^l}^{2^{l+1}} \(\Sigma_1(e^y; {\rm g}) +\Sigma_2(e^y; {\rm g})\)dy=o\( Y(\log Y)^{2m+2n-2} \).
\end{eqnarray*}
This completes the proof of Lemma \ref{pf-distri-Err-gene-lem}. \qed

\subsection{Lower order sum and error from the truncation}\label{subsec-lower-trunc}
In this section, we examine the lower order sum and the error from the truncation by a fixed large $T_0$, and give the proof of Lemma \ref{lem-distri-extra}. 

For the lower order sum, by (\ref{pf-distri-S-1}), we have
\begin{equation*}
\int_2^Y \left| S_1(e^y; \chi) \right|^2 dy\ll \sum_{j=2}^k (\log Y)^{2k-2j} \int_2^Y \left| \sum_{|\gamma|\leq T(e^y)}\frac{e^{i\gamma y}}{\frac{1}{2}+i\gamma} \right|^2 dy.
\end{equation*}
For the inner integral, by Lemma \ref{lem-Ford-S-2.2} and Lemma \ref{lem-Ford-S}, and the definition of $T=T(x)$,
\begin{align*}
\int_2^Y \left| \sum_{|\gamma|\leq T(e^y)}\frac{e^{i\gamma y}}{\frac{1}{2}+i\gamma} \right|^2 dy&\leq \sum_{l\leq \frac{\log Y}{\log 2}+1} \int_{2^l}^{2^{l+1}}\( \sum_{\substack{|\gamma_1-\gamma_2|\leq 1\\ |\gamma_1|, |\gamma_2|\leq T_{l'}}} +\sum_{\substack{|\gamma_1-\gamma_2|> 1\\ |\gamma_1|, |\gamma_2|\leq T_{l'}}} \) \frac{e^{i(\gamma_1-\gamma_2)y}}{(\frac{1}{2}+i\gamma_1)(\frac{1}{2}-i\gamma_2)} dy \nonumber\\
&\ll \sum_{l\leq \frac{\log Y}{\log 2}+1}  \(2^l\sum_{\gamma} \frac{\log |\gamma|}{|\gamma|^2}+ \sum_{\gamma_1, \gamma_2} \frac{1}{|\gamma_1| |\gamma_2| |\gamma_1-\gamma_2|}\)\ll Y.
\end{align*}
Thus,
\begin{equation*}
\int_2^Y \left| S_1(e^y; \chi) \right|^2 dy\ll \sum_{j=2}^k Y (\log Y)^{2k-2j}=o(Y (\log Y)^{2k-2})).
\end{equation*}

Next, we examine $S_2(x, T_0; \chi)$.  For fixed $T_0$, let $X_0$ be the largest $x$ such that $T=T(x)\leq T_0$. Since $x\leq T(x)\leq 2 x^2$, $\log X_0 \asymp \log T_0$. By Lemma \ref{lem-Ford-S-2.2} and Lemma \ref{lem-Ford-S},
\begin{eqnarray*}
&&\int_2^Y | S_2(e^y, T_0; \chi)|^2 dy\leq \int_2^{\log X_0} \left| \sum_{|\gamma|\leq T_0} \frac{1}{|\gamma|} \right|^2 dy +\int_{\log X_0}^Y \left| \sum_{T_0\leq |\gamma|\leq T(e^y)} \frac{e^{i\gamma y}}{\frac{1}{2}+i\gamma} \right|^2 dy \nonumber\\
&& \ll \log^5 T_0+\sum_{\frac{\log\log X_0}{\log 2}\leq l\leq \frac{\log Y}{\log 2}+1} \int_{2^l}^{2^{l+1}} \( \sum_{\substack{|\gamma_1-\gamma_2|\leq 1\\ T_0\leq |\gamma_1|, |\gamma_2|\leq T_{l'}}} +\sum_{\substack{|\gamma_1-\gamma_2|> 1\\ T_0\leq |\gamma_1|, |\gamma_2|\leq T_{l'}}} \) \frac{e^{i(\gamma_1-\gamma_2)y}}{(\frac{1}{2}+i\gamma_1)(\frac{1}{2}-i\gamma_2)} dy \nonumber\\
&&\ll \log^5 T_0+ \sum_{\frac{\log\log X_0}{\log 2}\leq l\leq \frac{\log Y}{\log 2}+1}  \(2^l\sum_{|\gamma|\geq T_0} \frac{\log |\gamma|}{|\gamma|^2}+ \sum_{|\gamma_1|, |\gamma_2| \geq T_0} \frac{1}{|\gamma_1||\gamma_2||\gamma_1-\gamma_2|}\)\nonumber\\
&& \ll Y \frac{\log^2 T_0}{T_0} +\log Y \frac{\log^3 T_0}{T_0}+\log^5 T_0.
\end{eqnarray*}
This completes the proof of this lemma. \qed

\section{Asymptotic formulas for the logarithmic densities}
In this section, we give the proof of Theorem \ref{thm-density-asym}.

For large $q$, Fiorilli and Martin \cite{F-M} gave an asymptotic formula for $\delta_{\Omega_1}(q; a, b)$. Lamzouri \cite{Lamz} also derived such an asymptotic formula using another method. Here, we want to derive asymptotic formulas for $\delta_{\Omega_k}(q; a, b)$ and $\delta_{\omega_k}(q; a, b)$ for fixed $q$ and large $k$.

Let $a$ be a quadratic non-residue $\bmod q$ and $b$ be a quadratic residue $\bmod q$, and $(a, q)=(b, q)=1$. Letting  $\lambda_k=\frac{1}{2^{k-1}}$, similar to formula (2.10) of \cite{F-M}, we have, under the assumptions ${\rm ERH_q}$ and ${ \rm LI_q}$,
\begin{equation*}
\delta_{\Omega_k}(q; a, b)=\frac{1}{2}+ \frac{(-1)^k}{2\pi}\int_{-\infty}^{\infty} \frac{\sin(\lambda_k (N(q; a)-N(q; b))x}{x} \Phi_{q; a, b}(x) dx.
\end{equation*}
Noting that $N(q, a)-N(q, b)=-A(q)$,
\begin{equation}\label{pf-den-asym-formu}
\delta_{\Omega_k}(q; a, b)=\frac{1}{2}+\frac{(-1)^{k-1}}{2\pi} \int_{-\infty}^{\infty}\frac{\sin(\lambda_k A(q) x)}{x} \Phi_{q; a, b}(x) dx.
\end{equation}
For any $\epsilon>0$,
\begin{equation}\label{den-asym-div}
\int_{-\infty}^{\infty}\frac{\sin(\lambda_k A(q) x)}{x} \Phi_{q; a, b}(x) dx=\(\int_{-\infty}^{\frac{1}{\lambda_k^{\epsilon}}}+\int_{-\frac{1}{\lambda_k^{\epsilon}}}^{\frac{1}{\lambda_k^{\epsilon}}}+\int_{\frac{1}{\lambda_k^{\epsilon}}}^{\infty} \) \frac{\sin(\lambda_k A(q) x)}{x} \Phi_{q; a, b}(x) dx.
\end{equation}
By Proposition 2.17 in \cite{F-M}, $|\Phi_{q; a, b}(t)|\leq e^{-0.0454\phi(q)t}$ for $t\geq 200$. So for large enough $k$,
\begin{equation}\label{pf-den-asym-est-1}
\int_{\frac{1}{\lambda_k^{\epsilon}}}^{\infty}\frac{\sin(\lambda_k A(q) x)}{x} \Phi_{q; a, b}(x) dx\ll \lambda_k  \int_{\frac{1}{\lambda_k^{\epsilon}}}^{\infty} e^{-0.0454\phi(q)x} dx \ll_{q, J, \epsilon}\lambda_k^{J}, ~\text{for any}~J>0.
\end{equation}
The integral over $x\leq -\frac{1}{\lambda_k^{\epsilon}}$ is also bounded by $\lambda_k^{J}$.

By Lemma 2.22 in \cite{F-M}, for each nonnegative integer $K$ and real number $C>1$, we have, uniformly for $|z|\leq C$,
\begin{equation*}
\frac{\sin z}{z}=\sum_{j=0}^{K} (-1)^j \frac{z^{2j}}{(2j+1)!}+O_{C, K}\(|z|^{2K+2} \).
\end{equation*}
Thus, the second integral in (\ref{den-asym-div}) is equal to
\begin{eqnarray}\label{pf-den-asym-est-ma}
&&\lambda_k A(q)\int_{-\frac{1}{\lambda_k^{\epsilon}}}^{\frac{1}{\lambda_k^{\epsilon}}} \frac{\sin(\lambda_k A(q) x)}{\lambda_k A(q)x} \Phi_{q; a, b}(x) dx\nonumber\\
&&=\sum_{j=0}^K \lambda^{2j+1}_k\frac{(-1)^j  A(q)^{2j+1}}{(2j+1)!} \int_{-\frac{1}{\lambda_k^{\epsilon}}}^{\frac{1}{\lambda_k^{\epsilon}}} x^{2j} \Phi_{q; a, b}(x) dx+O_{q, K}\(\lambda_k^{2K+3-\epsilon} \)\nonumber\\
&&=\sum_{j=0}^K \lambda^{2j+1}_k\frac{(-1)^j  A(q)^{2j+1}}{(2j+1)!} \int_{-\infty}^{\infty} x^{2j} \Phi_{q; a, b}(x) dx+O_{q, K, \epsilon}\(\lambda_k^{2K+3-\epsilon} \).
\end{eqnarray}

Combining (\ref{pf-den-asym-formu}), (\ref{pf-den-asym-est-1}), and (\ref{pf-den-asym-est-ma}), we get the asymptotic formula (\ref{formu-delta-Omega}) for $\delta_{\Omega_k}(q; a, b)$. Similarly, or by the results in Theorem \ref{thm-density-exist}, we have the asymptotic formula (\ref{formu-delta-omega}) for $\delta_{\omega_k}(q; a, b)$. \qed

\section{The source of main terms and proof of Lemma \ref{lem-main-term}}\label{sec-main-lem}

In this section, we give the proof of the main lemma we used for extracting out the bias terms and oscillating terms from the integrals over Hankel contours. 

Let $\mathcal{H}(0, X)$ be the truncated Hankel contour surrounding $0$ with radius $r$. Lau and Wu \cite{Lau-Wu} proved the following lemma. 
\begin{lem}[\cite{Lau-Wu}, Lemma 5]\label{lem-Gamma-j} For $X>1$, $z\in \mathbb{C}$ and $j\in \mathbb{Z}^{+}$, we have
	\begin{equation*}
		\frac{1}{2\pi i}\int_{\mathcal{H}(0, X)} w^{-z}(\log w)^j e^w dw=(-1)^j \frac{d^j}{dz^j}\(\frac{1}{\Gamma(z)}\)  +E_{j, z}(X),
	\end{equation*}
	where
	\begin{equation*}
		|E_{j,z}(X)|\leq \frac{e^{\pi|\Im(z)|}}{2\pi}\int_{X}^{\infty} \frac{(\log t+\pi)^j}{t^{\Re(z)}e^t}dt.
	\end{equation*}
\end{lem}

\noindent\textit{\textbf{Proof of Lemma \ref{lem-main-term}.}}
 We have the equality
\begin{equation*}
\frac{1}{s}=\frac{1}{a}+\frac{a-s}{a^2}+\frac{(a-s)^2}{a^2 s}.
\end{equation*}
With the above equality, we write the integral in the lemma as
\begin{equation*}
\frac{1}{2\pi i}\int_{\mathcal{H}(a, \delta)} \log^k(s-a)\(\frac{1}{a}+\frac{a-s}{a^2}+\frac{(a-s)^2}{a^2 s}\)x^s ds =: I_1+I_2+I_3.
\end{equation*}
For $I_3$, using Lemma \ref{lem-Esigma}, we get
\begin{align}\label{pf-lem-main-I3}
&\int_{\mathcal{H}(a,\delta)} \log^k(s-a)\frac{(a-s)^2}{a^2 s} x^s ds\nonumber\\
&\leq \left|\int_r^{\delta}\left((\log \sigma-i\pi)^k-(\log\sigma+i\pi)^k\right)\sigma^2 x^{-\sigma}\frac{x^a}{a^2 (a-\sigma)}d\sigma\right|\nonumber\\
&\quad +\int_{-\pi}^{\pi} x^{\Re(a)+r} \(\log \frac{1}{r}+\pi\)^k \frac{r^2}{|a|^2 |\Re(a)-r|}  r d\alpha\nonumber\\
&\ll\frac{|x^a|}{|a|^2 |\Re(a)-\delta|} \(\int_0^{\delta} |(\log \sigma -i\pi)^k-(\log\sigma +i\pi)^k|\sigma^2 x^{-\sigma}d\sigma+ \frac{(\log\frac{1}{r}+\pi)^{k}}{(1/r)^3}\)\nonumber\\
&\ll_k \frac{|x^a|}{|a|^2|\Re(a)-\delta|}\(\frac{(\log\log x)^{k-1}}{(\log x)^3}+ \frac{1}{x^{3-\epsilon}}\)\ll_k \frac{|x^a|}{|a|^2|\Re(a)-\delta|}\frac{(\log\log x)^{k-1}}{(\log x)^3}.
\end{align}
We estimate $I_2$ similarly. By Lemma \ref{lem-Esigma},
\begin{eqnarray}\label{pf-lem-main-I2}
&&\int_{\mathcal{H}(a,\delta)} \log^k(s-a)\frac{a-s}{a^2} x^s ds\nonumber\\
&&\leq \left|\int_r^{\delta}\left((\log \sigma-i\pi)^k-(\log\sigma+i\pi)^k\right)\sigma x^{-\sigma}\frac{x^a}{a^2 }d\sigma\right|+\int_{-\pi}^{\pi} x^{\Re(a)+r} \(\log \frac{1}{r}+\pi\)^k \frac{r}{|a|^2 }  r d\alpha\nonumber\\
&&\ll\frac{|x^a|}{|a|^2 } \(\int_0^{\delta} |(\log \sigma -i\pi)^k-(\log\sigma +i\pi)^k|\sigma x^{-\sigma}d\sigma+ \frac{(\log\frac{1}{r}+\pi)^{k}}{(1/r)^2}\)\nonumber\\
&&\ll_k \frac{|x^a|}{|a|^2}\(\frac{(\log\log x)^{k-1}}{(\log x)^2}+ \frac{1}{x^{2-\epsilon}}\)\ll_k \frac{|x^a|}{|a|^2}\frac{(\log\log x)^{k-1}}{(\log x)^2}.
\end{eqnarray}
For $I_1$, using change of variable $(s-a)\log x=w$, by Lemma \ref{lem-Gamma-j}, we get
\begin{align}\label{pf-lem-main-I1}
I_1&=\frac{1}{2\pi i}\frac{1}{\log x}\int_{\mathcal{H}(0,\delta\log x)}(\log w-\log\log x)^k \frac{x^a e^w}{a}dw\nonumber\\
&=\frac{x^a}{a\log x}  (-1)^k (\log\log x)^k \frac{1}{2\pi i}\int_{\mathcal{H}(0,\delta\log x)} e^w dw \nonumber\\
& \quad +(-1)^{k-1}k \frac{x^a}{a\log x}(\log\log x)^{k-1} \frac{1}{2\pi i}\int_{\mathcal{H}(0, \delta\log x)}e^w\log w  dw \nonumber\\
& \quad +\frac{x^a}{a\log x}\sum_{j=2}^k {k\choose j} \frac{1}{2\pi i}  \int_{\mathcal{H}(0, \delta\log x)} (-\log\log x)^{k-j} (\log w)^j e^w dw  \nonumber\\
&=\frac{(-1)^k x^a}{a\log x} \left\{  k (\log\log x)^{k-1}+\sum_{j=2}^k {k \choose j} (\log\log x)^{k-j} \frac{1}{\Gamma_j(0)}  \right\}  \nonumber \\
&  \quad +\frac{x^a}{a\log x} \sum_{j=1}^k {k \choose j} E_{j,0}(\delta\log x) (-\log\log x)^{k-j}.
\end{align}
By Lemma \ref{lem-Gamma-j},
\begin{equation*}
	|E_{j,0}(\delta\log x)|\leq \frac{1}{2\pi}\int_{\delta\log x}^{\infty} \frac{(\log t+\pi)^j}{e^t}dt
	\ll_j e^{-\frac{\delta\log x}{2}} \int_{\frac{\delta\log x}{2}}^{\infty} \frac{(\log t)^j}{e^{t/2}}dt
	\ll_j x^{-\frac{\delta}{2}}.
\end{equation*}
Hence, we get
\begin{equation}\label{pf-lem-main-E}
\left| \frac{x^a}{a\log x} \sum_{j=1}^k {k \choose j} E_{j,0}(\delta\log x) (-\log\log x)^{k-j} \right| \ll_k \frac{x^{\Re(a)}}{|a|\log x} \sum_{j=1}^k x^{-\frac{\delta}{2}} (\log\log x)^{k-j}\ll_k \frac{|x^{a-\delta/3|}}{|a|}.
\end{equation}

Combining (\ref{pf-lem-main-I3}), (\ref{pf-lem-main-I2}), (\ref{pf-lem-main-I1}), and (\ref{pf-lem-main-E}), we get the conclusion of Lemma \ref{lem-main-term}. \qed


\vspace{1em}
\textbf{Acknowledgments.} This research is partially
supported by NSF grants DMS-1201442 and DMS-1501982. I would like to thank my advisor, Professor Kevin Ford, for his kindly encouragement, useful suggestions and financial support to finish this project. I am grateful for the helpful comments of Dr. Youness Lamzouri. I also thank the encouragement of my friend Junjun Cheng.  The author would like to thank the referee for helpful comments.

{\footnotesize
Department of Mathematics, University of Illinois at Urbana-Champaign, 1409 West Green Street, Urbana, IL 61801, USA

\textit{E-mail}:  xmeng13@illinois.edu,


\begin{thebibliography}{99}
\bibitem{Che}P. L. Chebyshev, Lettre de M. le professeur Tch{\'e}byshev \'a M. Fuss, sur un nouveau th{\'e}oreme r{\'e}latif aux nombres premiers contenus dans la formes $4n+1$ et $4n+3$, \textit{Bull. de la Classe phys.-math. de l'Acad. Imp. des Sciences St. Petersburg} \textbf{11} (1853), 208.
\bibitem{Dave}H. Davenport, \textit{Multiplicative number theory, 3rd ed.}, Graduate Texts in Mathematics, vol. \textbf{74}, Springer-Verlag, New York-Berlin, 2000.
\bibitem{Diri}L. Dirichlet, Beweis des Satzes, Da\ss jede unbegrenzte arithmetische Progression $\ldots$ unendlich viele Primzahlen enth{\"a}lt, \textit{Abh. K{\"o}nig. Preuss. Akad.,} \textbf{34} (1837), 45-81. Reprinted on pp. 313-342 in \textit{Dirichlets Werke,} vol. \textbf{1}, Reimer, Berlin, 1889-97 and Chelsea, Bronx (NY), 1969.
\bibitem{F-M}D. Fiorilli and G. Martin, Inequalities in the Shanks-R\'{e}nyi prime number race: An asymptotic formula for the densities, \textit{J. reine angew. Math.}, \textbf{676} (2013), 121-212.
\bibitem{F-K}K. Ford and S. Konyagin, Chebyshev's conjecture and the prime number race. \textit{IV international Conference "Modern Problems of Number Theory and its Applications": Current Problems,} Part II (Russian) (Tula, 2001), 67-91, Mosk. Gos. Univ. im. Lomonosova, Mekh-Mat. Fak., Moscow, 2002.
\bibitem{Ford-S}K. Ford, J. Sneed, Chebyshev's bias for products of two primes. \textit{Experiment. Math.}, Volume \textbf{19}, Issue \textbf{4} (2010), 385-398.
\bibitem{Gran-Mar}A. Granville and G. Martin, Prime number races, \textit{Amer. Math. Monthly} \textbf{113} (2006), No. \textbf{1}, 1-33.

\bibitem{Hudson}R. Hudson, A common combinatorial principle underlies Riemann's formula,
the Chebyshev phenomenon, and other subtle effects in comparative prime number
theory. I. , \textit{J. reine angew. Mat}. \textbf{313} (1980), 133-150.  







\bibitem{Kara}A. A. Karatsuba, \textit{Basic Analytic Number Theory,} Springer-Verlag, 1993.
\bibitem{KT}S. Knapowski and P Tur{\'a}n, Comparative Prime Number Theory I., \textit{Acta. Math. Sci. Hungar.} \textbf{13} (1962), 315-342.
\bibitem{Lamz}Y. Lamzouri, Prime number races with three or more competitors, \textit{Math. Ann.}, (2013) \textbf{356}: 1117-1162.
\bibitem{Lan}E. Landau, Handbuch der Lehre von der Verteilung der Primzahlen (2 vols.), Teubner, Leipzig; 3rd edition: Chelsea, New York (1974).
\bibitem{Lau-Wu}Y. K. Lau, J. Wu. Sums of some multiplicative functions over a special set of integers. \textit{Acta Arithmetica.} \textbf{101.4} (2002)
\bibitem{Lee}J. Leech, Note on the distribution of prime numbers, \textit{J. London Math. Soc.} \textbf{32} (1957), 56-58.
\bibitem{Lit}J. E. Littlewood, Sur la distribution des nombres premiers, \textit{C. R. Acad. des Sciences Paris} \textbf{158} (1914), 1869-1872.
\bibitem{Mac}I.G. Macdonald, \textit{Symmetric functions and Hall polynomials, 2nd ed.}, Oxford Mathematical Monographs, Oxford University Press, New York, 1995.
\bibitem{Me-Re}A. Mendes and J. Remmel. \textit{Counting with Symmetric Functions}, Developments in Mathematics, volume \textbf{43}, Cham: Springer, 2015.
\bibitem{RS}M. Rubinstein and P. Sarnak, Chebyshev's bias, \textit{Experiment. Math.} \textbf{3} (1994), 173-197.
\bibitem{Tene}G. Tenenbaum, \textit{Introduction to analytic and probabilistic number theory, 3rd ed.,} Graduate studies in mathematics, vol. \textbf{163}, Providence, Rhode Island: American Mathematical Society, 2015.
\bibitem{Titch}E. C. Titchmarsh, \emph{The Theory of the Riemann Zeta-function, 2nd ed. rev. by D. R. Heath-Brown}, Clarendon Press, Oxford 1986




\end{thebibliography}
\end{document}